\documentclass[leqno,11pt]{amsart}
\usepackage{amsfonts,amsmath,amssymb,graphicx,psfrag}

\newtheorem{theorem}[equation]{Theorem}
\newtheorem{corollary}[equation]{Corollary}
\newtheorem{lemma}[equation]{Lemma}
\newtheorem{proposition}[equation]{Proposition}
\newtheorem{remark}[equation]{Remark}

\newtheorem{example}{Example}[section]
\numberwithin{equation}{section}

\begin{document}

\newcommand\1{\text{\rm 1}\hspace{-.95mm}\text{\rm I}}
\newcommand\C{\mathbb{C}}
\renewcommand\epsilon{\varepsilon}
\newcommand\N{\mathbb{N}}
\newcommand\R{\mathbb{R}}
\newcommand\ssf{\hspace{.25mm}}
\newcommand\ssb{\hspace{-.25mm}}
\newcommand\supp{\operatorname{supp}}
\newcommand\Z{\mathbb{Z}}
\newcommand\const{\operatorname{const.}}

\title{The Hardy space $H^1$ in the rational Dunkl setting}

\author[J.-Ph. Anker]{Jean-Philippe Anker}
\address{Universit\'e d'Orl\'eans \& CNRS,
F\'ed\'eration Denis Poisson (FR 2964) \& Laboratoire MAPMO (UMR 6628),
B\^atiment de math\'ematiques,
B.P. 6759, 45067 Orl\'eans cedex 2, France}
\email{anker@univ-orleans.fr}

\author[N. Ben Salem]{N\'ejib Ben Salem}
\address{Universit\'e de Tunis El Manar,
Facult\'e des Sciences de Tunis,
LR11ES11 Analyse Math\'ematiques et Applications,
2092, Tunis, Tunisie}
\email{Nejib.Bensalem@fst.rnu.tn}

\author[J. Dziuba\'nski]{Jacek Dziuba\'nski}
\address{Uniwersytet Wroc\l awski,
Instytut Matematyczny,
Pl. Grunwaldzki 2/4,
50-384 Wroc\l aw,
Poland}
\email{Jacek.Dziubanski@math.uni.wroc.pl}

\author[N. Hamda]{Nabila Hamda}
\address{Universit\'e de Tunis El Manar,
Facult\'e des Sciences de Tunis,
LR11ES11 Analyse Math\'ematiques et Applications,
2092, Tunis, Tunisie}
\email{hamdanabila@gmail.com}



\thanks{JPA and JD partially supported by the European Commission
(European program ToK \textit{HANAP\/} 2005-2009,
Polish doctoral program SSDNM 2010-2013)
and by a French-Polish cooperation program
(PHC Polonium 22529RF 2010-2011).
JPA, NBS, NHK partially supported
by a French-Tunisian cooperation program
(PHC Utique\,/\,CMCU 10G1503 2010-2013).
JD partially supported by the Polish National Science Centre
(Narodowe Centrum Nauki, grant DEC-2012/05/B/ST1/00672)
and by the University of Orl\'eans.}

\date{\today}

\begin{abstract}
This paper consists in a first study of the Hardy space $H^1$ in the rational Dunkl setting.
Following Uchiyama's approach, we characterizee $H^1$
atomically and by means of the heat maximal operator.
We also obtain a Fourier multiplier theorem for $H^1$.
These results are proved here in the one-dimensional case and in the product case.
\end{abstract}

\maketitle

\section{Introduction}

Dunkl theory is a far reaching generalization of Euclidean Fourier analysis,
which includes most special functions related to root systems,
such as spherical functions on Riemannian symmetric spaces.
It started in the late eighties with Dunkl's seminal article \cite{Dunkl}
and developed extensively afterwards.
We refer to the lecture notes \cite{Roesler3} for the rational Dunkl theory,
to the lecture notes \cite{Opdam} for the trigonometric Dunkl theory,
and to the books \cite{Cherednik, Macdonald} for the generalized quantum theories.

This paper deals with the real Hardy space $H^1$ in the rational Dunkl setting,
where the underlying space is of homogeneous type in the sense of Coifman-Weiss.
In such a setting, the theory of Hardy spaces goes back to the seventies
\cite{CoifmanWeiss, MaciasSegovia}.
Here we follow Uchyama's approach \cite{Uchiyama}
and we characterize the Hardy space $H^1$ in two ways,
by means of the heat maximal operator and atomically.
The first characterization,
which requires precise heat kernel estimates,
has lead us to a seemingly new observation,
namely that the heat kernel has a rather slow decay in certain directions
and is in particular not Gaussian in the present setting
(see Remark \ref{DecayRemark}).
The second characterization is used to prove a Fourier multiplier theorem for $H^1$.

Throughout the paper we shall restrict to the one-dimensional case and to the product case.
This restriction is due to our present lack of knowledge in general
about the behavior of the Dunkl kernel on the one hand
and about generalized translations on the other hand.

After this informal introduction,
let us introduce some notation and state our main results.
On \ssf$\R^n$ we consider the Dunkl operators
\begin{equation*}
D_jf(\mathbf{x})
=\tfrac{\partial}{\partial x_j} f(\mathbf{x})
+\tfrac{k_j}{x_j}\bigl[f(\mathbf{x})\ssb-\ssb f(\sigma_j\ssf\mathbf{x})\bigr]
\qquad(j\!=\!1,2,\dots,n)
\end{equation*}
associated with the reflections
\begin{equation}\label{Reflections}
\sigma_j\ssf(x_1,x_2,\dots,x_j,\dots,x_n)=(x_1,x_2,\dots,-x_j,\dots,x_n)
\end{equation}
and the multiplicities \ssf$k_j\!\ge\ssb0$\ssf.
Their joint eigenfunctions constitute the Dunkl kernel
\begin{equation}\label{DunklKernelProduct}
\mathbf{E}(\mathbf{x},\mathbf{y})
=\prod\nolimits_{\ssf j=1}^{\,n}\ssb E_{k_j}(x_j,y_j)\,,
\end{equation}
where
\begin{equation}\begin{aligned}\label{DunklKernel1D}
E_k(x,y)&=\tfrac{\Gamma(k\ssf+\frac12)}{\Gamma(k)\,\Gamma(\frac12)}
\int_{-1}^{+1}\hspace{-1mm}du\,(1\!-\ssb u)^{k-1}\ssf(1\!+\ssb u)^k\,e^{\,x\ssf y\ssf u}\\
&=e^{\,x\ssf y}
\underbrace{
\tfrac{\Gamma(2\ssf k\ssf+1)}{\Gamma(k)\,\Gamma(k\ssf+1)}\ssf
\int_{\ssf0}^{\,1}\!dv\hspace{.7mm}
v^{\ssf k-1}\ssf(1\!-\ssb v)^k\,e^{-2\ssf x\ssf y\ssf v}
}_{\textstyle{}_1F_{\ssf1}(k\ssf;2\ssf k\ssb+\!1\ssf;-\ssf2\,x\ssf y)}
\end{aligned}\end{equation}
(see for instance \cite[Example\;2.34]{Roesler3}).
Here \,${}_1F_{\ssf1}(a\ssf;b\ssf;z)$ \,is the confluent hypergeometric function,
which is also known as the Kummer function and denoted by $M(a,b,z)$\ssf.
Notice that \,$\mathbf{E}(\mathbf{x},\mathbf{y})\ssb
=\ssb e^{\ssf\langle\mathbf{x},\mathbf{y}\rangle}$
\ssf if all multiplicities \ssf$k_j$ vanish.

Let us first define the Hardy space $H^1$ by means of the heat maximal operator.
The Dunkl laplacian
\begin{equation*}
\mathbf{L}f(\mathbf{x})
=\sum\nolimits_{\ssf j=1}^{\,n}\ssb D_j^{\ssf2}f(\mathbf{x})
=\sum\nolimits_{\ssf j=1}^{\,n}\ssf\Bigl\{
\bigl(\tfrac\partial{\partial\ssf x_j}\bigr)^2\ssb f(\mathbf{x})\ssb
+\ssb\tfrac{2\ssf k_j}{x_j}\ssf\tfrac{\partial}{\partial\ssf x_j}f(\mathbf{x})\ssb
-\ssb\tfrac{k_j}{x_j^2}\ssf\bigl[\ssf f(\mathbf{x})\!
-\!f(\sigma_j\ssf\mathbf{x})\ssf\bigr]\Bigr\}
\end{equation*}
is the infinitesimal generator of the heat semigroup
\begin{equation*}
e^{\,t\,\mathbf{L}}\qquad(\ssf t\!>\!0\ssf)\ssf,
\end{equation*}
which acts by linear self-adjoint operators on $L^2(\R^n,d\boldsymbol{\mu})$
and by linear contractions on $L^p(\R^n,d\boldsymbol{\mu})\ssf$,
for every $1\ssb\le\ssb p\ssb\le\ssb\infty$\ssf,
where
\begin{equation}\label{ProductMeasure}
d\boldsymbol{\mu}(\mathbf{x})
=d\mu_1(x_1)\ssf\dots\,d\mu_n(x_n)
=|x_1|^{\ssf2\ssf k_1}\dots\,|x_n|^{\ssf2\ssf k_n}\,
dx_1\ssf\dots\ssf dx_n
\end{equation}
The heat semigroup consists of integral operators
\begin{equation*}
e^{\,t\ssf\mathbf{L}}f(\mathbf{x})
=\int_{\ssf\R^n}\!d\boldsymbol{\mu}(\mathbf{y})\,
\mathbf{h}_{\ssf t}(\mathbf{x},\mathbf{y})\,f(\mathbf{y})\,
\end{equation*}
associated with the heat kernel \cite{Roesler2}
\begin{equation}\label{HeatKernelProduct}
\mathbf{h}_{\ssf t}(\mathbf{x},\mathbf{y})=
\mathbf{c}_{\ssf\mathbf{k}}^{-1}\,t^{-\frac{\mathbf{N}}2}\,
e^{-\frac{|\mathbf{x}|^2+\ssf|\mathbf{y}|^2}{4\,t}}\,
\mathbf{E}\bigl(\tfrac{\mathbf{x}}{\sqrt{2\ssf t\ssf}},
\tfrac{\mathbf{y}}{\sqrt{2\ssf t\ssf}}\bigr)\ssf,
\end{equation}
where
\begin{equation}\label{HomogeneousDimension}
\mathbf{N}=n+\sum\nolimits_{\ssf j=1}^{\,n}\ssb2\,k_j
\end{equation}
is the homogeneous dimension and
\begin{equation*}
\mathbf{c}_{\ssf\mathbf{k}}\ssf=\,2^{\frac{\mathbf{N}}2}\ssb
\int_{\ssf\R^n}\!d\boldsymbol{\mu}(\mathbf{x})\,e^{-\frac{|\mathbf{x}|^2}2}\ssf
=\,2^{\ssf\mathbf{N}}\,\prod\nolimits_{\ssf j=1}^{\,n}\ssb\Gamma(k_j\!+\!\tfrac12)\,.
\end{equation*}
From this point of view,
the Hardy space \ssf$H^1$ consists of all functions
\ssf$f\!\in\!L^1(\R^n,d\boldsymbol{\mu})$
\ssf whose maximal heat transform
\begin{equation}\label{MaximalHeatOperator}
\mathbf{h}_{\ssf*}f(\mathbf{x})
=\ssf\sup\nolimits_{\,t>0}\,\Bigl|\ssf\int_{\ssf\R^n}\!d\boldsymbol{\mu}(\mathbf{y})\,
\mathbf{h}_{\ssf t}(\mathbf{x},\mathbf{y})\,f(\mathbf{y})\ssf\Bigr|
\end{equation}
belongs to \ssf$L^1(\R^n,d\boldsymbol{\mu})$
and the norm is given by
\begin{equation*}
\|f\|_{H^1}=\ssf\|\ssf\mathbf{h}_{\ssf*}f\ssf\|_{L^1}\ssf.
\end{equation*}

Let us turn next to the atomic definition of the Hardy space $H^1$.
Notice that $\R^n$,
equipped with the Euclidean distance
\ssf$d\ssf(\mathbf{x},\mathbf{y})\ssb=\ssb|\ssf\mathbf{x}\!-\!\mathbf{y}\ssf|$
\ssf and with the measure \ssf$\boldsymbol{\mu}$\ssf,
is a space of homogeneous type in the sense of Coifman-Weiss (see Appendix A).
Recall that an atom is a measurable function \,$a\ssb:\ssb\R^n\!\to\ssb\C$
\,such that
\begin{itemize}
\item[$\bullet$]
$\;a$ \ssf is supported in a ball \ssf$B$\ssf,
$\vphantom{\displaystyle\int}$

\item[$\bullet$]
$\;\|a\|_{L^\infty}\!\lesssim\boldsymbol{\mu}\ssf(B)^{-1}$\ssf,

\item[$\bullet$]
$\;\displaystyle\int_{\ssf\R^n}\hspace{-1mm}d\boldsymbol{\mu}(\mathbf{x})\,a(\mathbf{x})=0$\ssf.
\end{itemize}
By definition, the atomic Hardy space \ssf$H^1_{\text{atom}}$
\ssf consists of all functions \ssf$f\!\in\!L^1(\R^n,d\boldsymbol{\mu})$
which can be written as
\ssf$f\hspace{-.4mm}=\ssb\sum_{\ssf\ell}\ssb\lambda_{\ssf\ell}\,a_{\ssf\ell}$\ssf,
where the \ssf$a_\ell$'s are atoms and
\ssf$\sum_{\ssf\ell}\ssb|\lambda_{\ssf\ell}|\!<\!+\infty$\ssf,
and the norm is given by
\begin{equation*}
\|f\|_{H^1_\text{atom}}=\,\inf\,\sum\nolimits_{\ssf\ell}|\lambda_\ell|\,,
\end{equation*}
where the infimum is taken over all atomic decompositions of \ssf$f$.

Our first main result is the following theorem.

\begin{theorem}\label{Theorem1}
The spaces \ssf$H^1$ \ssb and \ssf$H^1_{\text{\rm atom}}$ coincide
and their norms are equivalent, i.e., there exists a constant \,$C\!>\!0$ such that
\begin{equation*}
C^{-1}\,\|f\|_{H^1}\ssb\le\|f\|_{H^1_{\text{\rm atom}}}\!\le C\,\| f\|_{H^1}\ssf.
\end{equation*}
\end{theorem}

The Fourier transform in the Dunkl setting is given by
\begin{equation}\label{FourierTransform}
\mathcal{F}\ssb f(\boldsymbol{\xi})
=\ssf\mathbf{c}_{\ssf\mathbf{k}}^{-1}\!
\int_{\ssf\R^n}\!d\boldsymbol{\mu}(\mathbf{x})\,f(\mathbf{x})\,
\mathbf{E}(\mathbf{x},\ssb-\ssf i\ssf\boldsymbol{\xi})\,.
\end{equation}
It is an isometric isomorphism of \ssf$L^2(\R^n,d\boldsymbol{\mu})$ onto itself
and the inversion formula reads
\begin{equation*}
f(\mathbf{x})=\mathcal{F}^{\ssf2\ssb}f(-\mathbf{x})\,.
\end{equation*}
Notice that, if all multiplicities \ssf$k_j$ vanish,
then \eqref{FourierTransform} boils down to the classical Fourier transform
\begin{equation*}
\widehat{f}(\boldsymbol{\xi})
=\ssf (2\ssf\pi)^{-\frac n2}\!\int_{\ssf\R^n}\!d\mathbf{x}\,f(\mathbf{x})\,
e^{-\ssf i\ssf\langle\ssf\mathbf{x},\ssf\boldsymbol{\xi}\ssf\rangle}\,.
\end{equation*}
Our second main result is the following H\"ormander type multiplier theorem
(see \cite{Hormander} for the original multiplier theorem on $L^p$ spaces).

\begin{theorem}\label{Theorem2}
Let \,$\chi\ssb=\ssb\chi(\boldsymbol{\xi})$
be a smooth radial function on \,$\R^n$ such that
\begin{equation*}
\chi(\boldsymbol{\xi})=\begin{cases}
\,1&\text{if \;}|\boldsymbol{\xi}|\!\in\!\bigl[\frac12,2\ssf\bigr]\ssf,\\
\,0&\text {if \;}|\boldsymbol{\xi}|\!\notin\!\bigl(\frac14,4\ssf\bigr)\ssf.\\
\end{cases}\end{equation*}
If a function \;$m\ssb=\ssb m(\boldsymbol{\xi})$ \ssf on \,$\R^n$ satisfies
\begin{equation*}
M=\,\sup\nolimits_{\,t>0}\,
\|\,\chi\,m(t\,.\,)\ssf\|_{\hspace{.1mm}W_{\ssf2}^{\ssf\mathbf{N}/2\ssf+\ssf\epsilon}}
<+\infty\,,
\end{equation*}
for some \,$\epsilon\!>\!0$\ssf,
then the multiplier operator
\begin{equation*}
\mathcal{T}_{\ssf m\ssf}f
=\mathcal{F}^{-1}\{\ssf m\,(\mathcal{F}\ssb f)\}
\end{equation*}
is bounded on the Hardy space \ssf$H^1$ and
\begin{equation*}
\|\,\mathcal{T}_{\ssf m}\,\|_{H^1\to\ssf H^1}\lesssim\,M\ssf.
\end{equation*}
\end{theorem}

Here \ssf$W_{\ssf2}^{\ssf\sigma\vphantom{|}}(\R^n)$
denotes the classical $L^2$ Sobolev space on \ssf$\R^n$,
whose norm is given by
\begin{equation*}
\|\ssf g\ssf\|_{\hspace{.1mm}W_{\ssf2}^{\ssf\sigma\vphantom{|}}}
=\Bigl\{\,\int_{\ssf\R^n}\hspace{-1mm}d\mathbf{x}\,
(1\!+\!|\mathbf{x}|^2)^\sigma\,|\ssf\widehat{g}(\mathbf{x})|^2\,\Bigr\}^{1/2}\,.
\end{equation*}
Notice that the multiplier \ssf$m$ \ssf is continuous and bounded,
as \,$\frac{\mathbf{N}}2\ssb+\ssb\epsilon\ssb>\ssb\frac n2$\ssf.

The theory of classical real Hardy spaces in $\mathbb R^n$ originates
from the study of holomorphic functions of one variable in the upper half-plane.
We refer the reader to the original works
of Stein-Weiss \cite{SteinWeiss},
Burkholder-Gundy-Silverstein \cite{BurkholderGundySilverstein}
and Fefferman-Stein \cite{FeffermanStein}.
An important contribution to this theory lies
in the atomic decomposition introduced by Coifman \cite{Coifman}
and extended to spaces of homogeneous type
by Coifman-Weiss \cite{CoifmanWeiss} (see also \cite{MaciasSegovia}).
More information can be found in the book \cite{Stein}
and references therein.

Our paper is organized as follows.
Section \ref{HeatKernelEstimates1D} is devoted to the heat kernel in dimension $1$.
There we analyze its behavior thoroughly and we remove a small part,
in order to get Gaussian estimates similar to the Euclidean setting.
These results are extended to the product case
in Section \ref{HeatKernelEstimatesProduct}.
Section 4 is devoted to the proof of Theorem \ref{Theorem1}
and Section 5 to the proof of Theorem \ref{Theorem2}.
Section 6 consists of 3 appendices.
Appendix A contains information about the measure of balls,
which is used throughout the paper.
Appendices B and C are devoted to so-called folklore results
in connection with Uchiyama's Theorem,
which have been used for instance in
\cite{DziubanskiPreisnerWrobel}.

This paper results from two independent research works,
which were carried out by the first and third authors,
respectively by the second and fourth authors,
and which have been merged into a joint article.
\medskip

\section{Heat kernel estimates in dimension $1$}
\label{HeatKernelEstimates1D}

Consider first the one-dimensional Dunkl kernel
\ssf$E(x,y)\ssb=\ssb E_k(x,y)$\ssf.
As the case \ssf$k\ssb=\ssb0$ \ssf is trivial,
we may assume that \ssf$k\ssb>\ssb0$\ssf.

\begin{lemma}\label{PropertiesDunklKernel}
\begin{itemize}
\item[(a)]
$\ssf E(x,y)$ is a holomorphic function of \ssf$(x,y)\!\in\!\C^2$.
\item[(b)]
$\ssf E(x,y)\ssb>\ssb0$ \,for every \,$x,y\!\in\!\R$\ssf.
\item[(c)]
$\ssf E(x,y)$ has the following symmetry and rescaling properties\,{\rm:}
\begin{equation*}
\begin{cases}
\,E(x,y)=E(y,x)
\qquad\forall\;x,y\!\in\!\C\ssf,\\
\,E(\lambda\ssf x,y)=E(x,\lambda\ssf y)
\qquad\forall\;\lambda,x,y\!\in\!\C\ssf.\\
\end{cases}
\end{equation*}
\item[(d)]
For every \,$y\!\in\!\C$\ssf,
\,$x\ssb\mapsto\ssb E(x,y)$ is an eigenfunction of the Dunkl operator
\begin{equation*}
Df(x)=f^{\ssf\prime}(x)+\tfrac kx\,\{\ssf f(x)\!-\!f(-x)\ssf\}
\end{equation*}
and of the Dunkl laplacian
\begin{equation*}
Lf(x)=D^{\ssf2\ssb}f(x)
=f^{\ssf\prime\prime}(x)
+\tfrac{2\ssf k}x\ssf f^{\ssf\prime}(x)
-\tfrac k{x^2}\,\{\ssf f(x)\!-\!f(-x)\ssf\}\,.
\end{equation*}
More precisely
\begin{equation*}
D_xE(x,y)=y\,E(x,y)
\quad\text{and}\quad
L_xE(x,y)=y^{\ssf2}\ssf E(x,y)\,.
\end{equation*}
\item[(e)]
As \,$x\ssf y\ssb\to\ssb0$\ssf,
\begin{equation*}
E(x,y)=1+\text{\rm O}\ssf(|\ssf x\ssf y\ssf|)\,.
\end{equation*}
\item[(f)]
As \,$x\ssf y\ssb\to\ssb+\infty$\ssf,
\begin{equation*}
E(x,y)=\tfrac{2^{\ssf k}\,\Gamma(k\ssf+\frac12)}{\sqrt{\pi\ssf}}\,
e^{\ssf x\ssf y}\,(x\ssf y)^{-k}\,
\bigl\{\ssf1\ssb-\ssb\tfrac{k^2}2\ssf\tfrac1{x\ssf y}\ssb
+\text{\rm O}\ssf(\tfrac1{x^2y^2})\ssf\bigr\}\,.
\end{equation*}
\item[(g)]
As \,$x\ssf y\ssb\to\ssb-\infty$\ssf,
\begin{equation*}
E(x,y)=\tfrac{2^{\ssf k-1}\ssf k\,\Gamma(k\ssf+\frac12)}{\sqrt{\pi\ssf}}\,
e^{-\ssf x\ssf y}\,(-\ssf x\ssf y)^{-k-1}\,
\bigl\{\ssf1\ssb+\ssb\tfrac{k^2-1}2\ssf\tfrac1{x\ssf y}\ssb
+\text{\rm O}\ssf(\tfrac1{x^2y^2})\ssf\bigr\}\,.
\end{equation*}
\end{itemize}
\end{lemma}

\begin{proof}
The first four properties are known to hold in general.
In dimension $1$, they can be also deduced
from the explicit expression \eqref{DunklKernel1D},
as does (e).
As already observed in \cite[Section\;2]{RoeslerDeJeu}
(see also \cite[Example\;5.1]{Roesler3}),
the asymptotics of $E(x,y)$ at infinity follow from
the asymptotics of the confluent hypergeometric function,
which read, let say for \ssf$0\ssb<\ssb a\ssb<\ssb b$\ssf,
\begin{equation*}
{}_1F_{\ssf1}(a\ssf;b\ssf;z)\ssf
\sim\ssf\tfrac{\Gamma(b)}{\Gamma(a)}\;e^{\,z}\,z^{\ssf a-b}\,
\sum\nolimits_{\ssf\ell=0}^{+\infty}\ssb
\tfrac{(1-\ssf a)_\ell\ssf(b\ssf-\ssf a)_\ell}{\ell\,!}\,z^{-\ell}
\end{equation*}
as \ssf$z\ssb\to\ssb+\infty$ \ssf and
\begin{equation*}
{}_1F_{\ssf1}(a\ssf;b\ssf;z)\ssf\sim\ssf
\tfrac{\Gamma(b)}{\Gamma(b\ssf-\ssf a)}\,|z|^{-a}\,
\sum\nolimits_{\ssf\ell=0}^{+\infty}\ssb
\tfrac{(a)_\ell\ssf(a\ssf-\ssf b\ssf+1)_\ell}{\ell\,!}\,|z|^{-\ell}
\end{equation*}
as \ssf$z\ssb\to\ssb-\infty$ \ssf
(see for instance \cite[(13.5.1)]{AbramowitzStegun} or \cite[(13.7.2)]{NIST}).
\end{proof}

Consider next the one-dimensional heat kernel
\begin{equation}\label{HeatKernel1D}
h_{\ssf t}(x,y)
=c_{\ssf k}^{-1}\,t^{-k-\frac12}\,e^{-\frac{x^2+\ssf y^2}{4\ssf t}}\ssf
E(\tfrac x{\sqrt{2\ssf t\ssf}},\tfrac y{\sqrt{2\ssf t\ssf}})
=c_{\ssf k}^{-1}\,t^{-k-\frac12}\,e^{-\frac{(x-y)^2}{4\ssf t}}\ssf
{}_1F_{\ssf1}(k\ssf;2\ssf k\ssb+\!1\ssf;-\tfrac{x\ssf y}t)\ssf,
\end{equation}
where \,$c_{\ssf k}\hspace{-.4mm}=\ssb2^{\ssf2\ssf k+1}\hspace{.5mm}\Gamma(k\!+\!\frac12)$\ssf.

\begin{proposition}\label{PropertiesHeatKernel1D}
\begin{itemize}
\item[(a)]
$\ssf h_{\ssf t}(x,y)$ is a \,$C^\infty$ function
of \,$(t,x,y)\!\in\!(0,+\infty)\!\times\ssb\R^2$.
\item[(b)]
$\ssf h_{\ssf t}(x,y)\ssb>\ssb0$ \,for every \,$t\!>\!0$ and \,$x,y\!\in\!\R$\ssf.
\item[(c)]
$\ssf h_{\ssf t}(x,y)$ has the following symmetry and rescaling properties\,{\rm:}
\begin{equation*}
\begin{cases}
\;h_{\ssf t}(x,y)=h_{\ssf t}(y,x)
\qquad\forall\;x,y\!\in\!\R\ssf,\\
\;h_{\ssf\lambda^2t\ssf}(\lambda\ssf x,\lambda\ssf y)
=|\lambda|^{-2\ssf k-1}\,h_{\ssf t}(x,y)
\qquad\forall\;\lambda\!\in\!\R^*\ssb,\,\forall\;t\!>\!0\ssf,\,\forall\;x,y\!\in\!\R\ssf.\\
\end{cases}
\end{equation*}
\item[(d)]
$\ssf h_{\ssf t}(x,y)$ satisfies the heat equation
\begin{equation*}
\begin{cases}
\;\partial_{\ssf t}\ssf h_{\ssf t}(x,y)=L_y\ssf h_{\ssf t}(x,y)\ssf,\\
\;\lim_{\,t\searrow0}\ssf h_{\ssf t}(x,y)\,|y|^{\ssf2\ssf k}\ssf dy=\delta_x(y)\ssf.\\
\end{cases}
\end{equation*}
\item[(e)]
The heat kernel has the following global behavior\,{\rm:}
\begin{equation*}
h_{\ssf t}(x,y)\,\asymp\,\begin{cases}
\;t^{-k-\frac12}\,
e^{-\frac{x^2\ssb+\ssf y^2}{4\ssf t}}
&\text{if \;}|\ssf x\ssf y\ssf|\!\le\ssb t\ssf,\\
\;t^{-\frac12}\,(x\ssf y)^{-k}\,
e^{-\frac{(x-y)^2}{4\ssf t}}
&\text{if \;}x\ssf y\ssb\ge\ssb t\ssf,\\
\;t^{\ssf\frac12}\,(-\ssf x\ssf y)^{-k-1}\,
e^{-\frac{(x+y)^2}{4\ssf t}}
&\text{if \,}-\ssb x\ssf y\ssb\ge\ssb t\ssf,\\
\end{cases}
\end{equation*}
and the following asymptotics\,{\rm:}
\begin{equation*}
h_{\ssf t}(x,y)\,=\,\begin{cases}
\;c_{\ssf k}^{-1}\,t^{-k-\frac12}\,
e^{-\frac{x^2\ssb+\ssf y^2}{4\ssf t}}\,
\bigl\{\ssf1\ssb+\text{\rm O}\bigl(\frac{|x\ssf y|}t\bigr)\bigr\}
&\text{if \;}\frac{x\ssf y}t\!\to\ssb0\ssf,\\
\;\frac1{2\ssf\sqrt{\pi\ssf}}\;e^{-\frac{(x-y)^2}{4\ssf t}}\,
t^{-\frac12}\,(x\ssf y)^{-k}\,
\bigl\{\ssf1\ssb-k^2\ssf\frac t{x\ssf y}\ssb
+\text{\rm O}\ssf(\frac{t^2}{x^2y^2})\bigr\}
&\text{if \;}\frac{x\ssf y}t\!\to\!+\infty\ssf,\\
\;\frac k{2\ssf\sqrt{\pi\ssf}}\;e^{-\frac{(x+y)^2}{4\ssf t}}\,
t^{\ssf\frac12}\,(-\ssf x\ssf y)^{-k-1}\,
\bigl\{\ssf1\ssb+\text{\rm O}\ssf(-\frac t{x\ssf y})\bigr\}
&\text{if \;}\frac{x\ssf y}t\!\to\!-\infty\ssf.
\end{cases}
\end{equation*}
\item[(f)]
The following gradient estimates hold for the heat kernel\,{\rm:}
\begin{equation*}
\bigl|\tfrac\partial{\partial y}\ssf h_{\ssf t}(x,y)\bigr|\,
\lesssim\,\begin{cases}
\;t^{-k-\frac32}\ssf(\ssf|x|\!+\!|y|\ssf)\,
e^{-\frac{x^2\ssb+\ssf y^2}{4\ssf t}}
&\text{if \;}|\ssf x\ssf y\ssf|\!\le\ssb t\ssf,\\
\,\bigl\{\ssf t^{-\frac32}\ssf|\ssf x\!-\!y\ssf|\hspace{-.5mm}
+\ssb t^{-\frac12}\ssf|y|^{-1}\bigr\}\,
(x\ssf y)^{-k}\,e^{-\frac{(x-y)^2}{4\ssf t}}
&\text{if \;}x\ssf y\ssb\ge\ssb t\ssf,\\
\,\bigl\{\ssf t^{-\frac12}\ssf|\ssf x\!+\!y\ssf|\hspace{-.5mm}
+\ssb t^{\ssf\frac12}\ssf(\ssf|x|^{-1}\!+\ssb|y|^{-1})\bigr\}\,
(-\ssf x\ssf y)^{-k-1}\,
e^{-\frac{(x+y)^2}{4\ssf t}}
&\text{if \,}-\ssb x\ssf y\ssb\ge\ssb t\ssf.\\
\end{cases}
\end{equation*}
\end{itemize}
\end{proposition}

\begin{proof}
The first five properties follow from the expression \eqref{HeatKernel1D}
and from Lemma \ref{PropertiesDunklKernel}.
Let us turn to the proof of (f).
By differentiating \eqref{HeatKernel1D} with respect to \ssf$y$ \ssf
and by using the well-known formula
\begin{equation*}
\tfrac d{dz}\,{}_1F_{\ssf1}(a\ssf;b\ssf;z)
=\tfrac ab\,{}_1F_{\ssf1}(a\ssb+\!1\ssf;b\ssb+\!1\ssf;z)
\end{equation*}
(see for instance \cite[(13.4.8)]{AbramowitzStegun} or \cite[(13.3.15)]{NIST}),
we get
\vspace{-1mm}
\begin{equation*}
\tfrac\partial{\partial y}\ssf h_{\ssf t}(x,y)
=c_{\ssf k}^{-1}\,t^{-k-\frac12}\,e^{-\frac{(x\ssf-\ssf y)^2}{4\ssf t}}\ssf
\bigl\{\ssf\tfrac{x-y}{2\,t}\,{}_1F_{\ssf1}(k\ssf;2\ssf k\ssb+\!1\ssf;-\tfrac{x\ssf y}t)
-\tfrac k{2\ssf k\ssf+1}\ssf\tfrac xt\,
{}_1F_{\ssf1}(k\ssb+\!1\ssf;2\ssf k\ssb+\ssb2\ssf;-\tfrac{x\ssf y}t)\bigr\}\ssf.
\end{equation*}
We conclude by using again the behavior of the confluent hypergeometric function.
\end{proof}

\begin{remark}\label{DecayRemark}
It follows from Proposition 2.3.(e) and Appendix A that
\begin{equation*}
h_{\ssf t}(x,x)\ssf\asymp\ssf
\mu\ssf(B(x,\ssb\sqrt{t\,}))^{-1}
\quad\text{and}\quad
\,h_{\ssf t}(x,\ssb-\ssf x)\ssf
\asymp\ssf\mu\ssf(B(x,\ssb\sqrt{t\,}))^{-1}\ssf\tfrac t{t\,+\,x^2}
\end{equation*}
for every \,$t\!>\!0$ and \,$x\!\in\!\R$\ssf.
Observe in particular that the heat kernel has no global Gaussian behavior
and decays rather slowly in certain directions.
This phenomenon is even more striking
in the product case \eqref{HeatKernelProductFormula},
where
\begin{equation*}
\mathbf{h}_{\ssf t}(\mathbf{x},{\mathbf{y}})\ssf
\asymp\ssf\boldsymbol{\mu}(B(\mathbf{x},\sqrt{t\,}))^{-1}\ssf
\tfrac t{t\;+\,|\ssf\mathbf{x}\ssf-\ssf\mathbf{y}\ssf|^2}
\end{equation*}
if \,$t\!>\!0$\ssf, $\mathbf{x}\!\in\!\R^n$
and \,$\mathbf{y}\ssb=\ssb(-\ssf x_1,x_2,\ssf\dots,x_n)$\ssf.
\end{remark}

Let us eventually introduce a variant of the heat kernel with a Gaussian behavior.
Given two smooth bump functions
\ssf$\chi_1$ and \ssf$\chi_2$ on \ssf$\R$
\ssf such that
\begin{equation*}
\begin{cases}
\,0\ssb\le\ssb\chi_1\!\le\!1\ssf,\\
\,\chi_1\!=\!1\text{ \,on \,}\bigl[-1,+\frac12\ssf\bigr]\ssf,\\
\,\supp\chi_1\!\subset\!\bigl[-2\ssf,+\frac23\ssf\bigr]\ssf,
\end{cases}
\quad\text{and}\qquad
\begin{cases}
\,0\ssb\le\ssb\chi_2\!\le\!1\ssf,\\
\,\chi_2\!=\!1\text{ \,on \,}\bigl[\ssf0,+\frac12\ssf\bigr]\ssf,\\
\,\supp\chi_2\ssb\subset\!\bigl[-1,+1\ssf\bigr]\ssf,
\end{cases}
\end{equation*}
consider the smooth cutoff function
\begin{equation*}
\chi_t(x,y)=\begin{cases}\,
\chi_1\bigl(\tfrac{x+y}{x}\bigr)\,\chi_2\bigl(\tfrac t{x^2}\bigr)
&\text{if \;}x\!\ne\!0\,,\\
\hspace{13mm}0
&\text{if \;}x\!=\!0\,,\\
\end{cases}\end{equation*}
and the truncated heat kernel
\begin{equation*}
H_t(x,y)=\{\ssf1\ssb-\chi_t(x,y)\}\,h_{\ssf t}(x,y)
\qquad\forall\;t\!>\!0\ssf,\,\forall\;x,y\!\in\!\R\ssf.
\end{equation*}

\begin{remark}
The truncated heat kernel \,$H_t(x,y)$
inherits the following properties of
the heat kernel \,$h_{\ssf t}(x,y)$\,{\rm:}
\begin{itemize}
\item[(a)]
Smoothness\,{\rm:}
\ssf$H_{\ssf t}(x,y)$ is a \,$C^\infty$ function
of \,$(t,x,y)\!\in\!(0,+\infty)\!\times\ssb\R^2$.
\item[(b)]
Non-negativity\,{\rm:}
\ssf$H_{\ssf t}(x,y)\ssb\ge\ssb0$
\,for every \,$t\!>\!0$ \ssf and \,$x,y\!\in\!\R$\ssf.
\item[(c)]
Rescaling\,{\rm:}
\ssf$H_{\lambda^2t}(\lambda\ssf x,\lambda\ssf y)\ssb
=\ssb|\lambda|^{-2\ssf k-1}\ssf H_t(x,y)$
\,for every \,$\lambda\!\in\!\R^*$\ssb,
\ssf$t\!>\!0$ \ssf and \,$x,y\!\in\!\R$\ssf.
\item[(d)]
Approximation of identity\,{\rm:}
\,$\lim_{\,t\searrow0}\ssf H_{\ssf t}(x,y)\,|y|^{\ssf2\ssf k}\ssf dy=\delta_x(y)$
\,for every \,$x,y\!\in\!\R$\ssf.
\end{itemize}
\end{remark}

\begin{theorem}\label{EstimatesTruncatedHeatKernel1D}
The following estimates hold for the truncated heat kernel \,$H_t(x,y)$\ssf.
\begin{itemize}
\item[(a)]
On-diagonal estimate\,{\rm:}
\begin{equation*}
H_t(x,x)\asymp\mu\ssf(B(x,\ssb\sqrt{t\,}))^{-1}
\qquad\forall\;t\!>\!0\ssf,\,\forall\;x\!\in\!\R\ssf.
\end{equation*}
\item[(b)]
Off-diagonal Gaussian estimate\,{\rm:}
\vspace{-1mm}
\begin{equation*}
0\le H_t(x,y)\lesssim
\mu\ssf(B(x,\ssb\sqrt{t\,}))^{-1}\,
e^{-\frac{(x-y)^2}{c\,t}}
\qquad\forall\;t\!>\!0\ssf,\,\forall\;x,y\!\in\!\R\ssf.
\end{equation*}
\item[(c)]
Gradient estimate\,{\rm:}
\vspace{-1mm}
\begin{equation*}\textstyle
\bigl|\frac\partial{\partial y}H_t(x,y)\bigr|\lesssim
t^{-\frac12}\,\mu\ssf(B(x,\ssb\sqrt{t\,}))^{-1}\,
e^{-\frac{(x-y)^2}{c\,t}}
\qquad\forall\;t\!>\!0\ssf,\,\forall\;x,y\!\in\!\R\ssf.
\end{equation*}
\item[(d)]
Lipschitz estimates\,{\rm:}
\begin{equation*}\textstyle
|\ssf H_t(x,y)\ssb-\ssb H_t(x,y^{\ssf\prime})\ssf|
\lesssim\ssf\mu\ssf(B(x,\ssb\sqrt{t\,}))^{-1}\,
\frac{|\ssf y\ssf-\ssf y^{\ssf\prime}|}{\sqrt{t\,}}
\qquad\forall\;t\!>\!0\ssf,\,\forall\;x,y,y^{\ssf\prime}\hspace{-1mm}\in\!\R\ssf,
\end{equation*}
with the following improvement,
if \,$|\ssf y\!-\!y^{\ssf\prime}|\ssb\le\ssb\frac12\,|\ssf x\!-\!y\ssf|$\,{\rm:}
\vspace{-1mm}
\begin{equation*}\textstyle
|\ssf H_t(x,y)\ssb-\ssb H_t(x,y^{\ssf\prime})\ssf|
\lesssim\ssf\mu\ssf(B(x,\ssb\sqrt{t\,}))^{-1}\,
e^{-\frac{(x-y)^2}{c\,t}}\,
\frac{|\ssf y\ssf-\ssf y^{\ssf\prime}|}{\sqrt{t\,}}\ssf.
\end{equation*}
\end{itemize}
Here \,$c$ \ssf denotes some positive constant
and the ball measure has the following behavior,
according to Appendix A\,{\rm:}
\begin{equation*}
\mu\ssf(B(x,\ssb\sqrt{t\,}))\ssf
\asymp\ssf\begin{cases}
\;t^{\,k+\frac12}
&\text{if \,}|x|\ssb\le\!\sqrt{\ssf t\,},\\
\,|x|^{\ssf2\ssf k}\ssf\sqrt{\ssf t\,}
&\text{if \,}|x|\ssb\ge\!\sqrt{\ssf t\,}.\\
\end{cases}\end{equation*}
\end{theorem}

\begin{proof}
As far as (a), (b), (c) are concerned,
the case \ssf$x\ssb=\ssb0$ \ssf follows immediately
from the previous heat kernel estimates.
Thus we may assume that \ssf$x\ssb\ne\ssb0$
\ssf and reduce furthermore to \ssf$x\ssb=\ssb1$ \ssf by rescaling.
\medskip

\noindent
(a) is immediate\,:
\begin{equation*}
H_t(1,1)=h_{\ssf t}(1,1)
\asymp\left\{\begin{matrix}
t^{-\frac12}
&\text{if \,$t\ssb\le\!1$}\,\\
\,t^{-k-\frac12}
&\text{if \,$t\ssb\ge\!1$}\,\\
\end{matrix}\right\}\asymp\mu\ssf(B(1,\ssb\sqrt{t\,}))^{-1}\ssf.
\end{equation*}

\begin{figure}[ht]
\begin{center}
\psfrag{y}[c]{$y$}
\psfrag{t}[c]{$t$}
\psfrag{0}[c]{$0$}
\psfrag{1/2}[c]{$\frac12$}
\psfrag{-1/2}[c]{$-\frac12$}
\psfrag{-1}[c]{$-1$}
\psfrag{2}[c]{$2$}
\psfrag{-2}[c]{$-2$}
\includegraphics[height=80mm]{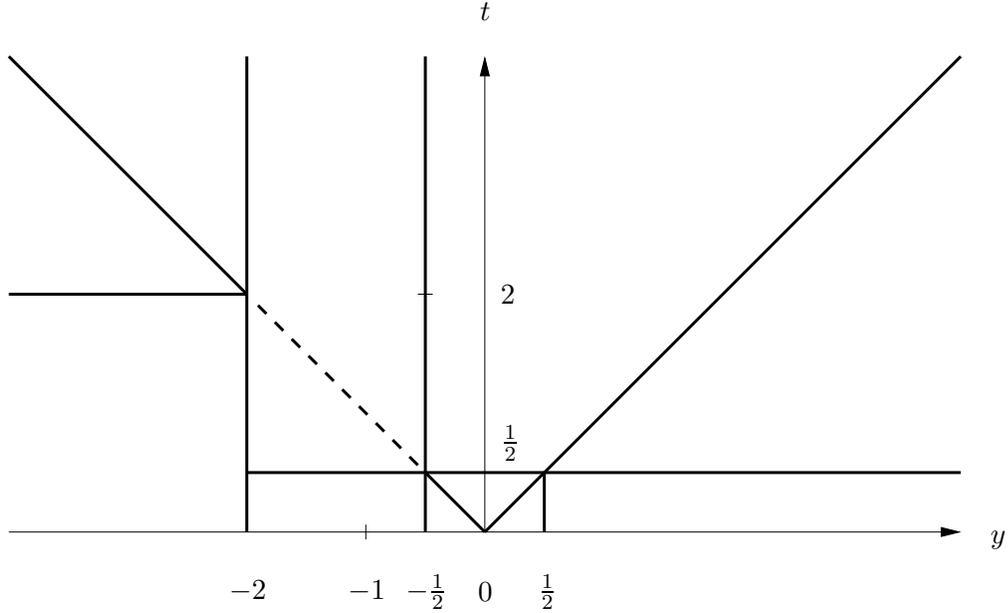}
\caption{Cases and subcases considered in the proofs of (b) and (c)}
\end{center}
\end{figure}

\noindent
Let us next prove (b).
\vspace{1mm}

\noindent$\bullet$
\,\textbf{Case 1.} Assume that \ssf$|y|\!\le\ssb t$\ssf.
\vspace{1mm}

\noindent$\circ$
\,\textbf{Subcase 1.1.}
Assume that $t$ \ssf is bounded above,
say \ssf$t\ssb\le\!\frac12$\ssf.
Then
\begin{equation*}
H_t(1,y)\le h_{\ssf t}(1,y)\asymp
t^{-k-\frac12}\,e^{-\frac{1+y^2}{4\ssf t}}\!
=t^{-\frac12}\,e^{-\frac{(1-y)^2}{8\ssf t}}\,
t^{-k}\,e^{-\frac{1+y^2}{8\ssf t}}\,e^{-\frac y{4\ssf t}}
\end{equation*}
is bounded above by
\begin{equation*}
\mu\ssf(B(1,\ssb\sqrt{t\,}))^{-1}
\,e^{-\frac{(1-y)^2}{8\hspace{.5mm}t}}
\end{equation*}
as \,$t^{\ssf\frac12}\ssb
\asymp\ssb\mu\ssf(B(1,\ssb\sqrt{t\,}))$\ssf,
\,$t^{-k}\ssb\lesssim e^{\ssf\frac1{8\ssf t}}\ssb
\le e^{\ssf\frac{1+y^2}{8\ssf t}}$
\ssf and \,$e^{\ssf\frac y{4\ssf t}}\ssb\asymp\ssb1$\ssf.
\vspace{1mm}

\noindent$\circ$
\,\textbf{Subcase 1.2.}
Assume that $t$ \ssf is bounded below,
say \ssf$t\ssb\ge\!\frac12$\ssf.
Then
\begin{equation*}
H_t(1,y)
\le h_{\ssf t}(1,y)
\asymp t^{-k-\frac12}\,e^{-\frac{1+y^2}{4\ssf t}}\!
=t^{-k-\frac12}\,e^{-\frac{(1-y)^2}{4\ssf t}}\,e^{-\frac y{2\ssf t}}
\end{equation*}
with \,$t^{\ssf k+\frac12}\ssb
\asymp\ssb\mu\ssf(B(1,\ssb\sqrt{t\,}))$
\ssf and \,$e^{\frac y{2\ssf t}}\ssb\asymp\ssb1$\ssf.
\vspace{1mm}

\noindent$\bullet$
\,\textbf{Case 2.}
Assume that $y$ is close to \ssf$-x\ssb=\!-1$\ssf,
say $-2\ssb\le\ssb y\ssb\le\!-\frac12$\ssf.

\vspace{1mm}

\noindent$\circ$
\,\textbf{Subcase 2.1.}
If \ssf$t\ssb\le\!\frac12\ssf(\ssf\le\!-\ssf y\ssf)$\ssf, then
\begin{equation*}
H_t(1,y)\ssb=0\,.
\end{equation*}

\noindent$\circ$
\,\textbf{Subcase 2.2.}
If $t$ \ssf is bounded below,
say \ssf$t\ssb\ge\!\frac12$\ssf,
we argue as in Subcase 1.2.
\vspace{1mm}

\noindent$\bullet$
\,\textbf{Case 3}.
Assume that \ssf$y\ssb\ge\ssb t$\ssf.
\vspace{1mm}

\noindent$\circ$
\,\textbf{Subcase 3.1.}
Assume that $t$ is bounded below,
say $(\ssf y\!\ge)\ssf t\!\ge\!\frac12$\ssf.
Then
\begin{equation*}\textstyle
H_t(1,y)\le h_{\ssf t}(1,y)\asymp
t^{-\frac12}\,y^{-k}\,e^{-\frac{(1-y)^2}{4\ssf t}}\ssb
\le t^{-k-\frac12}\,e^{-\frac{(1-y)^2}{4\ssf t}}
\end{equation*}
with \,$t^{\ssf k+\frac12}\ssb\asymp\mu\ssf(B(1,\ssb\sqrt{t\,}))$\ssf.
\vspace{1mm}

\noindent$\circ$
\,\textbf{Subcase 3.2.}
Assume that \ssf$y\ssb\ge\!\frac12\!\ge\ssb t$\ssf.
Then
\begin{equation*}\textstyle
H_t(1,y)\le h_{\ssf t}(1,y)\asymp
t^{-\frac12}\,y^{-k}\,e^{-\frac{(1-y)^2}{4\ssf t}}\ssb
\lesssim t^{-\frac12}\,e^{-\frac{(1-y)^2}{4\ssf t}}
\end{equation*}
with \,$t^{\ssf\frac12}\ssb\asymp\mu\ssf(B(1,\ssb\sqrt{t\,}))$\ssf.
\vspace{1mm}

\noindent$\circ$
\,\textbf{Subcase 3.3.}
Assume that \ssf$t\ssb\le\ssb y\ssb\le\!\frac12$\ssf.
Then
\begin{equation*}\textstyle
H_t(1,y)\le h_{\ssf t}(1,y)\asymp
t^{-\frac12}\,y^{-k}\,e^{-\frac{(1-y)^2}{4\ssf t}}
\end{equation*}
is bounded above by
\begin{equation*}
\mu\ssf(B(1,\ssb\sqrt{t\,}))^{-1}
\,e^{-\frac{(1-y)^2}{8\hspace{.5mm}t}}
\end{equation*}
as \,$t^{\ssf\frac12}\ssb\asymp\ssb\mu\ssf(B(1,\ssb\sqrt{t\,}))$
\ssf and \,$y^{-k}\ssb\le\ssb t^{-k}\ssb
\lesssim e^{\ssf\frac1{32\,t}}\ssb
\le e^{\frac{(1-y)^2}{8\ssf t}}$.
\vspace{1mm}

\noindent$\bullet$
\,\textbf{Case 4}.
Assume that \ssf$y\ssb\le\!-\ssf t\,(<\!0)$
and that $y$ \ssf stays away from \ssf$-1$\ssf,
say \ssf$y\!\notin\!\bigl(-2\ssf, -\frac12\ssf\bigr)$.
Notice that \ssf$(1\!+\!y)^2\!\ge\!\frac{(1\ssf-\,y)^2}9$
\,if and only if
\,$y\!\notin\!\bigl(-2\ssf, -\frac12\ssf\bigr)$.
\vspace{1.5mm}

\noindent$\circ$
\,\textbf{Subcase 4.1.}
Assume that \ssf$2\ssb\le\ssb t\ssb\le\!-\ssf y$\ssf.
Then
\begin{equation*}
H_t(1,y)\leq h_{\ssf t}(1,y)\asymp
t^{\ssf\frac12}\,(-\ssf y)^{-k-1}\,e^{-\frac{(1+y)^2}{4\ssf t}}\!
\le t^{-k-\frac12}\,e^{-\frac{(1-y)^2}{36\,t}}
\end{equation*}
with \,$t^{\ssf k+\frac12}\ssb\asymp\mu\ssf(B(1,\ssb\sqrt{t\,}))$\ssf.
\vspace{1mm}

\noindent$\circ$
\,\textbf{Subcase 4.2.}
Assume that \ssf$t\ssb\le\ssb2\ssb\le\!-\ssf y$\ssf.
Then
\begin{equation*}
H_t(1,y)\leq h_{\ssf t}(1,y)\asymp
t^{\ssf\frac12}\,(-\ssf y)^{-k-1}\,e^{-\frac{(1+y)^2}{4\ssf t}}\!
\lesssim t^{-\frac12}\,e^{-\frac{(1-y)^2}{36\,t}}
\end{equation*}
with \,$t^{\ssf\frac12}\ssb\asymp\mu\ssf(B(1,\ssb\sqrt{t\,}))$\ssf.
\vspace{1mm}

\noindent$\circ$
\,\textbf{Subcase 4.3.}
Assume that \ssf$t\ssb\le\!-\ssf y\ssb\le\ssb\frac12$\ssf.
Then
\vspace{-1mm}
\begin{equation*}
H_t(1,y)\leq h_{\ssf t}(1,y)\asymp
t^{\ssf\frac12}\,(-\ssf y)^{-k-1}\,e^{-\frac{(1+y)^2}{4\ssf t}}\!
\le t^{-k-\frac12}\,e^{-\frac{(1+y)^2}{8\ssf t}}\,e^{-\frac{(1-y)^2}{72\,t}}
\end{equation*}
\vspace{-3mm}

\noindent
is bounded above by
\vspace{-2mm}
\begin{equation*}
\mu\ssf(B(1,\ssb\sqrt{t\,}))^{-1}\,e^{-\frac{(1-y)^2}{72\,t}}
\end{equation*}
as \,$t^{\ssf\frac12}\ssb\asymp\mu\ssf(B(1,\ssb\sqrt{t\,}))$
\ssf and \,$t^{-k}\ssb\lesssim e^{\ssf\frac1{32\,t}}\ssb
\le e^{\frac{(1+y)^2}{8\ssf t}}$.
\medskip

\noindent
The proof of (c) follows the same pattern.
To begin with, observe that the derivative
\vspace{1mm}
\begin{equation*}\textstyle
\frac\partial{\partial y}\bigl\{1\!-\ssb\underbrace{
\chi_1(1\!+\ssb y)\,\chi_2(t)
\vphantom{\big|}}_{\chi_{\ssf t}(1,\ssf y)}\bigr\}
=-\,\chi_1^{\ssf\prime}(1\!+\ssb y)\,\chi_2(t)
\end{equation*}
of the cut-off is bounded and vanishes unless
\,$y\ssb\in\!\bigl(-3\ssf,-2\ssf\bigr)\ssb\cup\ssb\bigl(-\frac12,0\ssf\bigr)$
\,and \,$t\ssb\le\!1$.
According to the subcases 1.1, 4.2 and 4.3 above, the contribution of
\,$\frac\partial{\partial y}\ssf
\{\hspace{.2mm}1\hspace{-.4mm}-\hspace{-.2mm}\chi_{\hspace{.1mm}t}(1,y)\}\,
h_{\ssf t}(1,y)$ \ssf to \ssf$\frac\partial{\partial y}H_t(1,y)$
is bounded by
\vspace{-1mm}
\begin{equation*}
\mu\ssf(B(1,\ssb\sqrt{t\,}))^{-1}\,
e^{-\frac{(1-y)^2}{c\,t}}
\le t^{-\frac12}\,\mu\ssf(B(1,\ssb\sqrt{t\,}))^{-1}\,
e^{-\frac{(1-y)^2}{c\,t}}
\end{equation*}
Thus it remains for us to estimate the contribution of
\,$\{\hspace{.2mm}1\hspace{-.4mm}-\hspace{-.2mm}\chi_{\hspace{.1mm}t}(1,y)\}
\,\frac\partial{\partial y}\ssf h_{\ssf t}(1,y)$\ssf.
\vspace{1mm}

\noindent$\bullet$
\,\textbf{Case 1.} Assume that \ssf$|y|\!\le\ssb t$\ssf.
\vspace{1mm}

\noindent$\circ$
\,\textbf{Subcase 1.1.}
Assume that \ssf$t\ssb\le\!\frac12$\ssf.
Then
\vspace{-4mm}
\begin{align*}
\{\hspace{.2mm}1\hspace{-.4mm}-\hspace{-.2mm}\chi_{\hspace{.1mm}t}(1,y)\}\,
\bigl|\tfrac\partial{\partial y}\ssf h_{\ssf t}(1,y)\bigr|
&\lesssim t^{-k-\frac32}\ssf(1\!+\!|y|)\,e^{-\frac{1+y^2}{4\ssf t}}\!
\lesssim\ssf\overbrace{
t^{-k-\frac12}\ssf e^{-\frac{1+y^2}{8\ssf t}}\ssf e^{-\frac y{4\ssf t}}
}^{\text{bounded}}\ssf
t^{-1}\ssf e^{-\frac{(1-y)^2}{8\ssf t}}\\
&\lesssim
t^{-\frac12}\ssf
\mu\ssf(B(1,\ssb\sqrt{t\,}))^{-1}\,
e^{-\frac{(1-y)^2}{8\ssf t}}
\end{align*}
\vspace{-1mm}

\noindent$\circ$
\,\textbf{Subcase 1.2.}
Assume that \ssf$t\ssb\ge\!\frac12$\ssf.
Then
\vspace{-4mm}
\begin{align*}
\{\hspace{.2mm}1\hspace{-.4mm}-\hspace{-.2mm}\chi_{\hspace{.1mm}t}(1,y)\}\,
\bigl|\tfrac\partial{\partial y}\ssf h_{\ssf t}(1,y)\bigr|
&\lesssim t^{-k-\frac32}\ssf(1\!+\hspace{-.4mm}|y|\hspace{.1mm})\,
e^{-\frac{1+y^2}{4\ssf t}}
\!\lesssim t^{-k-1}\ssf e^{-\frac{(1-y)^2}{8\ssf t}}
\overbrace{
\bigl(\tfrac{1+\ssf y^2}t\bigr)^{\ssb\frac12}\ssf
e^{-\frac{1+y^2}{8\ssf t}}\ssf e^{-\frac y{4\ssf t}}
}^{\text{bounded}}\\
&\lesssim
t^{-\frac12}\ssf
\mu\ssf(B(1,\ssb\sqrt{t\,}))^{-1}\,
e^{-\frac{(1-y)^2}{8\ssf t}}
\end{align*}
\vspace{-1mm}

\noindent$\bullet$
\,\textbf{Case 2.}
Assume that $-2\ssb\le\ssb y\ssb\le\!-\frac12$\ssf.
\vspace{1mm}

\noindent$\circ$
\,\textbf{Subcase 2.1.}
If \ssf$t\ssb\le\!\frac12\ssf(\ssf\le\!-\ssf y\ssf)$\ssf, then
\begin{equation*}
\{\hspace{.2mm}1\hspace{-.4mm}-\hspace{-.2mm}\chi_{\hspace{.1mm}t}(1,y)\}\,
\tfrac\partial{\partial y}\ssf h_{\ssf t}(1,y)=0\,.
\end{equation*}

\noindent$\circ$
\,\textbf{Subcase 2.2.}
If $t$ \ssf is bounded below,
say \ssf$t\ssb\ge\!\frac12$\ssf,
we argue as in Subcase 1.2.
\vspace{1mm}

\noindent$\bullet$
\,\textbf{Case 3}.
Assume that \ssf$y\ssb\ge\ssb t$\ssf.
\vspace{1mm}

\noindent$\circ$
\,\textbf{Subcase 3.1.}
Assume that $(\ssf y\!\ge)\ssf t\!\ge\!\frac12$\ssf.
Then
\begin{align*}
\{\hspace{.2mm}1\hspace{-.4mm}-\hspace{-.2mm}\chi_{\hspace{.1mm}t}(1,y)\}\,
\bigl|\tfrac\partial{\partial y}\ssf h_{\ssf t}(1,y)\bigr|
&\lesssim\bigl\{\ssf t^{-\frac32}\ssf|1\hspace{-1mm}-\!y\ssf|\ssb
+t^{-\frac12}\ssf\ssf y^{-1}\bigr\}\,y^{-k}\,e^{-\frac{(1-y)^2}{4\ssf t}}\\
&\lesssim\ssf t^{-k-1}\,e^{-\frac{(1-y)^2}{8\ssf t}}
\underbrace{\Bigl
\{1\ssb+\ssb\tfrac{|1-\ssf y\ssf|}{\sqrt{t\,}}\,e^{-\frac{(1-y)^2}{8\ssf t}}
\Bigr\}}_{\text{bounded}}\\
&\lesssim\ssf t^{-\frac12}\ssf\mu\ssf(B(1,\ssb\sqrt{t\,}))^{-1}\,
e^{-\frac{(1-y)^2}{8\ssf t}}
\end{align*}
\vspace{-1mm}

\noindent$\circ$
\,\textbf{Subcase 3.2.}
Assume that \ssf$y\ssb\ge\!\frac12\!\ge\ssb t$\ssf.
Then
\begin{align*}
\{\hspace{.2mm}1\hspace{-.4mm}-\hspace{-.2mm}\chi_{\hspace{.1mm}t}(1,y)\}\,
\bigl|\tfrac\partial{\partial y}\ssf h_{\ssf t}(1,y)\bigr|
&\lesssim\bigl\{\ssf t^{-\frac32}\ssf|1\hspace{-1mm}-\!y\ssf|\ssb
+t^{-\frac12}\ssf y^{-1}\bigr\}\,y^{-k}\,e^{-\frac{(1-y)^2}{4\ssf t}}\\
&\lesssim\ssf t^{-1}\,e^{-\frac{(1-y)^2}{8\ssf t}}
\underbrace{\Bigl\{\sqrt{t\,}\!
+\ssb\tfrac{|1-\ssf y|}{\sqrt{t\,}}\,e^{-\frac{(1-y)^2}{8\ssf t}}
\Bigr\}}_{\text{bounded}}\\
&\lesssim\ssf t^{-\frac12}\ssf
\mu\ssf(B(1,\ssb\sqrt{t\,}))^{-1}\,
e^{-\frac{(1-y)^2}{8\ssf t}}
\end{align*}
\vspace{-1mm}

\noindent$\circ$
\,\textbf{Subcase 3.3.}
Assume that \ssf$t\ssb\le\ssb y\ssb\le\!\frac12$\ssf.
Then
\begin{align*}
\{\hspace{.2mm}1\hspace{-.4mm}-\hspace{-.2mm}\chi_{\hspace{.1mm}t}(1,y)\}\,
\bigl|\tfrac\partial{\partial y}\ssf h_{\ssf t}(1,y)\bigr|
&\lesssim\bigl\{\ssf t^{-\frac32}\ssf|1\hspace{-1mm}-\!y\ssf|\ssb
+t^{-\frac12}\ssf y^{-1}\bigr\}\,y^{-k}\,e^{-\frac{(1-y)^2}{4\ssf t}}\\
&\lesssim\ssf t^{-1}\,e^{-\frac{(1-y)^2}{8\ssf t}}
\underbrace{\vphantom{\big|}
t^{-k-\frac12}\,e^{-\frac1{32\ssf t}}
}_{\text{bounded}}\\
&\lesssim\ssf t^{-\frac12}\ssf
\mu\ssf(B(1,\ssb\sqrt{t\,}))^{-1}\,
e^{-\frac{(1-y)^2}{8\ssf t}}
\end{align*}
\vspace{-1mm}

\noindent$\bullet$
\,\textbf{Case 4}.
Assume that \ssf$y\ssb\le\!-\ssf t\,(<\!0)$
and that \ssf$y\!\notin\!\bigl(-2\ssf, -\frac12\ssf\bigr)$.
Recall that \ssf$(1\!+\!y)^2\!\ge\!\frac{(1\ssf-\,y)^2}9$
\,if and only if
\,$y\!\notin\!\bigl(-2\ssf, -\frac12\ssf\bigr)$.
\vspace{1.5mm}

\noindent$\circ$
\,\textbf{Subcase 4.1.}
Assume that \ssf$2\ssb\le\ssb t\ssb\le\!-\ssf y$\ssf.
Then
\begin{align*}
\{\hspace{.2mm}1\hspace{-.4mm}-\hspace{-.2mm}\chi_{\hspace{.1mm}t}(1,y)\}\,
\bigl|\tfrac\partial{\partial y}\ssf h_{\ssf t}(1,y)\bigr|
&\lesssim\bigl\{\ssf t^{-\frac12}\ssf|1\hspace{-1mm}+\!y\ssf|\ssb
+t^{\ssf\frac12}\tfrac{1+\ssf|y|}{|y|}\bigr\}\,
|y|^{-k-1}\,e^{-\frac{(1+y)^2}{4\ssf t}}\\
&\lesssim\ssf t^{-k-1}\,e^{-\frac{(1+y)^2}{8\ssf t}}
\underbrace{\tfrac{|1+\ssf y\hspace{.1mm}|}{\sqrt{t\,}}\,
e^{-\frac{(1+y)^2}{8\ssf t}}}_{\text{bounded}}\\
&\lesssim\ssf t^{-\frac12}\ssf
\mu\ssf(B(1,\ssb\sqrt{t\,}))^{-1}\,
e^{-\frac{(1-y)^2}{72\,t}}
\end{align*}
\vspace{-1mm}

\noindent$\circ$
\,\textbf{Subcase 4.2.}
Assume that \ssf$t\ssb\le\ssb2\ssb\le\!-\ssf y$\ssf.
Then
\begin{align*}
\{\hspace{.2mm}1\hspace{-.4mm}-\hspace{-.2mm}\chi_{\hspace{.1mm}t}(1,y)\}\,
\bigl|\tfrac\partial{\partial y}\ssf h_{\ssf t}(1,y)\bigr|
&\lesssim\bigl\{\ssf t^{-\frac12}\ssf|1\hspace{-1mm}+\!y\ssf|\ssb
+t^{\ssf\frac12}\tfrac{1+\ssf|y|}{|y|}\bigr\}\,
|y|^{-k-1}\,e^{-\frac{(1+y)^2}{4\ssf t}}\\
&\lesssim\ssf t^{-1}\,e^{-\frac{(1+y)^2}{8\ssf t}}
\underbrace{\Bigl\{\tfrac{|1+\ssf y|}{\sqrt{t\,}}\,e^{-\frac{(1+y)^2}{8\ssf t}}\!
+\ssb\sqrt{t\ssf}\ssf\Bigr\}}_{\text{bounded}}\\
&\lesssim\ssf t^{-\frac12}\ssf
\mu\ssf(B(1,\ssb\sqrt{t\,}))^{-1}\,
e^{-\frac{(1-y)^2}{72\,t}}
\end{align*}
\vspace{-1mm}

\noindent$\circ$
\,\textbf{Subcase 4.3.}
Assume that \ssf$t\ssb\le\!-\ssf y\ssb\le\ssb\frac12$\ssf.
Then
\begin{align*}
\chi_{\hspace{.1mm}t}(1,y)\,
\bigl|\tfrac\partial{\partial y}\ssf h_{\ssf t}(1,y)\bigr|
&\lesssim\bigl\{\ssf t^{-\frac12}\ssf|1\hspace{-1mm}+\!y\ssf|\ssb
+t^{\ssf\frac12}\tfrac{1+\ssf|y|}{|y|}\bigr\}\,
|y|^{-k-1}\,e^{-\frac{(1+y)^2}{4\ssf t}}\\
&\lesssim\ssf t^{-1}\,e^{-\frac{(1+y)^2}{8\ssf t}}
\underbrace{\vphantom{\big{|}}
t^{-k-\frac12}\,e^{-\frac1{32\ssf t}}
}_{\text{bounded}}\\
&\lesssim\ssf t^{-\frac12}\ssf
\mu\ssf(B(1,\ssb\sqrt{t\,}))^{-1}\,
e^{-\frac{(1-y)^2}{72\,t}}
\end{align*}
\medskip

\noindent
Eventually, (d) is an immediate consequence of (c).
For every
\ssf$y^{\ssf\prime\prime}\hskip-1mm\in\![\,y,y^{\ssf\prime}\ssf]$\ssf,
we have indeed
\vspace{-1mm}
\begin{equation*}
e^{-\frac{(x-y^{\ssf\prime\prime})^2}{c\,t}}\le1\,.
\end{equation*}
Moreover, if
\ssf$|\ssf y\!-\!y^{\ssf\prime}|\!\le\!\frac12\ssf|\ssf x\!-\!y\ssf|$\ssf,
then
\ssf$|\ssf x\!-\!y^{\ssf\prime\prime}|\ssb
\ge\ssb|\ssf x\!-\!y\ssf|\ssb-\ssb|\ssf y\!-\!y^{\ssf\prime\prime}|\ssb
\ge\ssb|\ssf x\!-\!y\ssf|\ssb-\ssb|\ssf y\!-\!y^{\ssf\prime}|\ssb
\ge\ssb\frac12\ssf|\ssf x\!-\!y\ssf|$\ssf,
hence
\vspace{-1mm}
\begin{equation*}
e^{-\frac{(x-y^{\ssf\prime\prime})^2}{c\ssf t}}
\le\ssf e^{-\frac{(x-y)^2}{4\ssf c\ssf t}}\,.
\end{equation*}
\end{proof}

\noindent
\begin{remark}\label{RemarkEstimatesTruncatedHeatKernel1D}
Contrarily to \ssf$h_{\ssf t}(x,y)$\ssf,
\ssf$H_t(x,y)$ \ssf is not symmetric in the space variables \ssf$x,y$\ssf.
Nevertheless, according to the following result,
we may replace \ssf$\mu\ssf(B(x,\ssb\sqrt{t\,}))$
\ssf by \ssf$\mu\ssf(B(y,\ssb\sqrt{t\,}))$
\ssf in the estimates (b), (c)
and in the second estimate (d).
\end{remark}

\begin{lemma}
For every \ssf$\varepsilon\!>\!0$\ssf,
there exists \ssf$C\!>\!0$ \ssf such that
\begin{equation*}
\tfrac{\mu\ssf(B(x,\sqrt{t\ssf}))}{\mu\ssf(B(y,\sqrt{t\ssf}))}\ssf
\le\ssf C\,e^{\,\varepsilon\ssf\frac{(x-y)^2}t}
\qquad\forall\;x,y\!\in\!\R\ssf,
\,\forall\;t\!>\!0\ssf.
\end{equation*}
\end{lemma}
\begin{proof}
By rescaling (see Appendix A),
we can reduce to the case \ssf$t\ssb=\!1$\ssf.
The estimate
\begin{equation*}
\tfrac{\mu\ssf(B(x,1))}{\mu\ssf(B(y,1))}\ssf
\lesssim\,e^{\ssf\varepsilon\ssf(x-y)^2}
\end{equation*}
is obvious if \ssf$x$ \ssf and \ssf$y$ \ssf are bounded
or if \ssf$|x|/|y|$ \ssf is bounded from above.
In the remaining case,
let say when \ssf$|x|\!\ge\!1\!+\hspace{-.4mm}2\ssf|y|$\ssf,
we have \ssf$|x|\ssb\le\ssb|x\!-\!y|\ssb+\ssb|y|\ssb
\le\ssb|x\!-\!y|\ssb+\ssb\frac12\ssf|x|$\ssf,
hence \ssf$|x|\ssb\le\ssb2\,|x\!-\!y|$\ssf.
Furthermore,
as \ssf$|x\!-\!y|\ssb\ge\ssb|x|\ssb-\ssb|y|\ssb\ge\ssb1$\ssf,
we have \ssf$|x|\ssb\le\ssb2\,(x\!-\!y)^2$.
Thus
\begin{equation*}
\tfrac{\mu\ssf(B(x,1))}{\mu\ssf(B(y,1))}\ssf
\lesssim\ssf\mu\ssf(B(x,1))\ssf
\asymp\ssf(\ssf|x|\!+\!1)^{2k}\ssf
\lesssim\,e^{\,\frac\varepsilon2\ssf|x|}
\lesssim\,e^{\,\varepsilon\ssf(x-y)^2}\ssf.
\end{equation*}
\end{proof}

Next proposition, which will be used in the proof of Theorem \ref{Theorem1},
shows that the truncated heat kernel \ssf$H_t(x,y)$ captures
the main features of the heat kernel \ssf$h_{\ssf t}(x,y)$.

\begin{proposition}\label{MaximalOperatorError1D}
The maximal operator
\begin{equation*}
Q_*f(x)=\ssf\sup\nolimits_{\ssf t>0}\,
\Bigl|\ssf\int_{\ssf\R}\ssb d\mu(y)\,Q_t(x,y)\,f(y)\ssf\Bigr|\ssf,
\end{equation*}
associated with the error
\begin{equation*}
Q_t(x,y)
=h_{\ssf t}(x,y)\ssb-\ssb H_t(x,y)
=\chi_{\hspace{.1mm}t}(x,y)\,h_{\ssf t}(x,y)
\ge0\,,
\end{equation*}
is bounded from \ssf$L^1(\R,d\mu)$ into itself.
\end{proposition}

\begin{proof}
It suffices to check that
\begin{equation*}
\sup\nolimits_{\,y\in\R}\ssb\int_{\ssf\R}d\mu(x)\,
\sup\nolimits_{\,t>0}Q_t(x,y)\ssf<+\infty\,.
\end{equation*}
The case \ssf$y\ssb=\ssb0$ \ssf is trivial,
as \ssf$\chi_{\hspace{.1mm}t}(x,0)$
and hence \ssf$Q_t(x,0)$ vanish,
for every \ssf$t\!>\!0$ and \ssf$x\!\in\!\R$\ssf.
Consider next the case \ssf$y\!\in\!\R^*$\ssb,
which reduces to \ssf$y\ssb=\ssb1$ by rescaling.
Then \ssf$\chi_{\hspace{.1mm}t}(x,1)$ and \ssf$Q_t(x,1)$ vanish,
unless \,$t\!<\!9$ \ssf and $-3\!<\!x\!<\!-\frac13$\ssf,
and in this range (see Proposition \ref{PropertiesHeatKernel1D})
\begin{equation*}
h_{\ssf t}(x,1)\asymp\ssf t^{\frac12}\,e^{-\frac{(x+1)^2}{4\ssf t}}
\end{equation*}
is bounded. Hence
\begin{equation*}
\int_{\ssf\R}d\mu(x)\,\sup\nolimits_{\,t>0}Q_t(x,1)
\ssf\lesssim\ssb\int_{-3}^{-\frac13}\hspace{-1mm}dx\,
\sup\nolimits_{\,0<t<9}h_{\ssf t}(x,1)
\ssf<+\infty\,.
\end{equation*}
\end{proof}

\section{Heat kernel estimates in the product case}
\label{HeatKernelEstimatesProduct}

According to \eqref{HeatKernelProduct} and \eqref{DunklKernelProduct},
the heat kernel in the product case splits up into one-dimensional heat kernels\,:
\begin{equation}\label{HeatKernelProductFormula}
{\mathbf{h}}_{\ssf t}(\mathbf{x},{\mathbf{y}})
=\prod\nolimits_{\ssf j=1}^{\,n}h_{\,t}^{(j)}(x_j,y_j)\,.
\end{equation}
By expanding
\begin{equation*}
h_{\,t}^{(j)}(x_j,y_j)\ssf=\,\underbrace{
\{\ssf 1\!-\ssb\chi_{\ssf t}(x_j,y_j)\}\,h_{\,t}^{(j)}(x_j,y_j)
}_{H_{\,t}^{\ssf(j)\ssb}(x_j,\,y_j)}\,
+\,\underbrace{
\chi_{\ssf t}(x_j,y_j)\,h_{\,t}^{(j)}(x_j,y_j)
}_{Q_{\,t}^{(j)\ssb}(x_j,\,y_j)}\,,
\end{equation*}
\vspace{-3mm}

\noindent
we get
\begin{equation*}
{\mathbf{h}}_{\ssf t}(\mathbf{x},{\mathbf{y}})
={\mathbf{H}}_{\ssf t}(\mathbf{x},{\mathbf{y}})
+{\mathbf{P}}_{\ssb t}(\mathbf{x},{\mathbf{y}})\ssf.
\end{equation*}
Here
\begin{equation*}
{\mathbf{H}}_{\ssf t}(\mathbf{x},{\mathbf{y}})
=\prod\nolimits_{\ssf j=1}^{\,n}H_{\,t}^{(j)}(x_j,y_j)
\end{equation*}
and \ssf${\mathbf{P}}_{\ssb t}(\mathbf{x},{\mathbf{y}})$
is the sum of all possible products
\begin{equation*}
\widetilde{\mathbf{P}}_{\ssb t}(\mathbf{x},\ssf{\mathbf{y}})
=\prod\nolimits_{\ssf j=1}^{\,n}p_{\,t}^{(j)}(x_j,y_j)\,,
\end{equation*}
where each factor \ssf$p_{\,t}^{(j)}(x_j,y_j)$ is equal
to \ssf$H_{\,t}^{(j)}(x_j,y_j)$ or \ssf$Q_{\,t}^{(j)}(x_j,y_j)$\ssf,
and at least one factor \ssf$p_{\,t}^{(j)}(x_j,y_j)$
is equal to \ssf$Q_{\,t}^{(j)}(x_j,y_j)$\ssf.
Notice the rescaling property
\begin{equation*}
{\mathbf{h}}_{\lambda^2t}(\lambda\ssf\mathbf{x},\lambda\ssf{\mathbf{y}})
=|\lambda|^{-\mathbf{N}}\,{\mathbf{h}}_{\ssf t}(\mathbf{x},{\mathbf{y}})
\qquad\forall\;\lambda\!\in\!\R^*\ssb,\,\forall\;t\!>\!0\ssf,
\,\forall\;\mathbf{x},\mathbf{y}\!\in\!\R^n,
\end{equation*}
and similarly for the other product kernels.
The following estimates follow from the one-dimensional case
(see Theorem \ref{EstimatesTruncatedHeatKernel1D}
and Remark \ref{RemarkEstimatesTruncatedHeatKernel1D}).

\begin{theorem}\label{EstimatesTruncatedHeatKernelProduct}
\begin{itemize}
\item[(a)]
On-diagonal estimate\,{\rm:}
\begin{equation*}
{\mathbf{H}}_{\ssf t}(\mathbf{x},\mathbf{x})\asymp
\boldsymbol{\mu}\ssf(\ssf\mathbf{B}\ssf(\ssf\mathbf{x},\ssb\sqrt{t\,}))^{-1}
\qquad\forall\;t\!>\!0\ssf,\,\forall\;\mathbf{x}\!\in\!\R^n.
\end{equation*}
\item[(b)]
Off-diagonal Gaussian estimate\,{\rm:}
\vspace{-1.5mm}
\begin{equation*}
0\le{\mathbf{H}}_{\ssf t}(\mathbf{x},{\mathbf{y}})
\lesssim\ssf\boldsymbol\max\ssf
\bigl\{\ssf\mu\ssf
(\ssf\mathbf{B}\ssf(\ssf\mathbf{x},\ssb\sqrt{t\,})),
\ssf\mu\ssf(\ssf\mathbf{B}\ssf
(\ssf{\mathbf{y}},\ssb\sqrt{t\,}))\bigr\}^{-1}\,
e^{-\frac{|\ssf\mathbf{x}\ssf-\ssf{\mathbf{y}}\ssf|^2}{c\,t}}
\end{equation*}
for every \,$t\!>\!0$ and for every \,$\mathbf{x},\mathbf{y}\!\in\!\R^n$.
\item[(c)]
Gradient estimate\,{\rm:}
\vspace{-1.5mm}
\begin{equation*}
|\ssf\nabla_{\ssb\mathbf{y}\ssf}
{\mathbf{H}}_{\ssf t}(\mathbf{x},{\mathbf{y}})\ssf|
\lesssim t^{-\frac12}\,
\boldsymbol\max\ssf
\bigl\{\ssf\mu\ssf
(\ssf\mathbf{B}\ssf(\ssf\mathbf{x},\ssb\sqrt{t\,})),
\ssf\mu\ssf(\ssf\mathbf{B}\ssf
(\ssf{\mathbf{y}},\ssb\sqrt{t\,}))\bigr\}^{-1}\,
e^{-\frac{|\ssf\mathbf{x}\ssf-\ssf{\mathbf{y}}\ssf|^2}{c\,t}}
\end{equation*}
for every \,$t\!>\!0$ and \,$\mathbf{x},\mathbf{y}\!\in\!\R^n$.
\item[(d)]
Lipschitz estimates\,{\rm:}
\begin{equation*}
|\ssf{\mathbf{H}}_{\ssf t}(\mathbf{x},{\mathbf{y}})\ssb
-\ssb{\mathbf{H}}_{\ssf t}(\mathbf{x},{\mathbf{y}}^{\ssf\prime})\ssf|
\lesssim\ssf\boldsymbol{\mu}\ssf(\ssf\mathbf{B}\ssf
(\ssf\mathbf{x},\ssb\sqrt{t\,}))^{-1}\,
\tfrac{|\ssf{\mathbf{y}}\ssf-\ssf{\mathbf{y}}^{\ssf\boldsymbol\prime}|}{\sqrt{t\,}}\,,
\end{equation*}
for every \,$t\!>\!0$ and
\,$\mathbf{x},\mathbf{y},\mathbf{y}^{\ssf\boldsymbol\prime}\hspace{-1mm}\in\!\R^n$,
with the following improvement,
if \,$|\ssf{\mathbf{y}}\!-\!{\mathbf{y}}^{\ssf\boldsymbol\prime}|\ssb
\le\ssb\frac12\,|\ssf\mathbf{x}\!-\!{\mathbf{y}}\ssf|$\,{\rm:}
\vspace{-1.5mm}
\begin{equation*}
|\ssf{\mathbf{H}}_{\ssf t}(\mathbf{x},{\mathbf{y}})\ssb
-\ssb{\mathbf{H}}_{\ssf t}(\mathbf{x},{\mathbf{y}}^{\ssf\boldsymbol\prime})\ssf|
\lesssim\ssf\boldsymbol\max\ssf
\bigl\{\ssf\mu\ssf(\ssf\mathbf{B}\ssf(\ssf\mathbf{x},\ssb\sqrt{t\,})),
\ssf\mu\ssf(\ssf\mathbf{B}\ssf(\ssf{\mathbf{y}},\ssb\sqrt{t\,}))\bigr\}^{-1}\,
e^{-\frac{|\ssf\mathbf{x}\ssf-\ssf{\mathbf{y}}\ssf|^2}{c\,t}}\,
\tfrac{|\ssf{\mathbf{y}}\ssf-\ssf{\mathbf{y}}^{\ssf\boldsymbol\prime}|}{\sqrt{t\,}}\,.
\end{equation*}
\end{itemize}
\end{theorem}

Let us turn to
the analog of Proposition \ref{MaximalOperatorError1D}
in the product case.

\begin{proposition}\label{MaximalOperatorErrorProduct}
The maximal operator
\begin{equation*}
{\mathbf{P}}_{\ssb*}f(\mathbf{x})
=\sup\nolimits_{\ssf t>0}\,
\Bigl|\ssf\int_{\ssf\R^n}\!d{\boldsymbol{\mu}}({\mathbf{y}})\,
{\mathbf{P}}_t(\mathbf{x},{\mathbf{y}})\,f({\mathbf{y}})\,\Bigr|\,,
\end{equation*}
is bounded from \ssf$L^1(\R^n\ssb,{\boldsymbol{\mu}})$ into itself.
\end{proposition}

\begin{proof}
We will show again that
\begin{equation*}
\sup\nolimits_{\,{\mathbf{y}}\in\R^n}
\int_{\ssf\R^n}\!d{\boldsymbol{\mu}}(\mathbf{x})\,
\sup\nolimits_{\,t>0}{\mathbf{P}}_{\ssb t}(\mathbf{x},{\mathbf{y}})\,
<\,+\infty\,,
\end{equation*}
but the proof will be more involved in the product case than in the one-dimensional case.
Let us begin with some observations.
First of all, by using the symmetries
\begin{equation*}
H_{\,t}^{(j)}(x_j,y_j)=H_{\,t}^{(j)}(-x_j,-y_j)
\quad\text{and}\quad
Q_{\,t}^{(j)}(x_j,y_j)=Q_{\,t}^{(j)}(-x_j,-y_j)
\end{equation*}
\vspace{-3mm}

\noindent
and by interchanging variables,
we may reduce to products of the form
\vspace{1mm}
\begin{equation*}
\widetilde{\mathbf{P}}_{\ssb t}(\mathbf{x},\ssf{\mathbf{y}})
=\underbrace{Q_{\,t}^{(1)}(x_1,y_1)\ssf\dots\,
Q_{\,t}^{(n^{\ssf\prime})}(x_{n^{\ssf\prime}}\ssb,y_{n^{\ssf\prime}})
}_{{\mathbf{Q}}_{\ssf t}^{\ssf\prime}
(\mathbf{x}^{\ssf\prime}\ssb,\ssf{\mathbf{y}}^{\ssf\prime})}
\underbrace{
H_{\,t}^{(n^{\ssf\prime}\ssb+1)}
(x_{n^{\ssf\prime}\ssb+1},y_{n^{\ssf\prime}\ssb+1})
\ssf\dots\,H_{\,t}^{(n)}(x_n,y_n)
}_{{\mathbf{H}}_{\ssf t}^{\ssf\prime\prime}
(\mathbf{x}^{\ssf\prime\prime}\ssb,\ssf{\mathbf{y}}^{\ssf\prime\prime})}
\end{equation*}
\vspace{-2.5mm}

\noindent
where \ssf$1\ssb\le\ssb n^{\ssf\prime}\!\le\ssb n$ \ssf and
\ssf$0\ssb\le\ssb y_1\!\le{\dots}\le\ssb y_{n^{\ssf\prime}}$\ssf.
Next we may assume that,
for every \ssf$1\ssb\le\ssb j\ssb\le\ssb n^{\ssf\prime}$,
\begin{equation*}
y_j\!>\ssb0\,,\quad
-\ssf3\hspace{.4mm}y_j\!<\ssb x_j\!<\!-\ssf\tfrac13\ssf y_j
\quad\text{and}\quad
x_j^2\ssb>\ssb t\ssf,
\end{equation*}
because otherwise
\ssf$\chi_{\ssf t}(x_j,y_j)$ and hence \ssf$Q_{\,t}^{(j)}(x_j,y_j)$
vanish.
Eventually, by rescaling, we may reduce to the case \ssf$y_1\!=\!1$\ssf.
Consequently,
\ssf$t$ \ssf is bounded by
\ssf$x_1^{\ssf2}\!<\ssb9\,y_1^{\ssf2}\!=\ssb9$
\ssf and each factor \ssf$Q_{\,t}^{(j)}(x_j,y_j)$ \ssf is bounded by
\begin{equation*}
t^\frac12\,(-\ssf x_j\ssf y_j)^{-k_j-1}\,e^{\ssf-\frac{(x_j+\ssf y_j)^2}{4\,t}}\,
\1_{\,\left(\ssb-\ssf3\ssf y_j,\ssf-\ssf\frac13\ssf y_j\ssb\right)}(x_j)\,
\lesssim\,t^\frac12\,y_j^{-2\ssf k_j-2}\,
\1_{\,\left(\ssb-\ssf3\ssf y_j,\ssf-\ssf\frac13\ssf y_j\ssb\right)}(x_j)\,.
\end{equation*}
Thus, on the one hand, the integral
\begin{equation*}\begin{aligned}
{\mathbf{I}}^{\ssf\prime}({\mathbf{y}}^{\ssf\prime})\ssf
&=\int_{\ssf\R^{n^\prime}}\!
d{\boldsymbol{\mu}}^{\ssf\prime}(\mathbf{x}^{\ssf\prime})\,
\sup\nolimits_{\,t>0\ssf}t^{-\frac{n'}2}\ssf{\mathbf{Q}}_{\,t}^{\ssf\prime}
(\mathbf{x}^{\ssf\prime}\ssb,{\mathbf{y}}^{\ssf\prime})\\
&\lesssim\int_{-\ssf3}^{-\frac13}\hspace{-1mm}d\mu_1(x_1)\,
y_2^{-2\ssf k_2-2}\!
\int_{-\ssf3\ssf y_2}^{-\frac13\ssf y_2}
\hspace{-1.5mm}d\mu_2(x_2)\;\dots\;
y_{n^{\ssf\prime}}^{-2\ssf k_{n^\prime}-2}\!
\int_{-\ssf3\ssf y_{n^{\ssf\prime}}}^{-\frac13\ssf y_{n^{\ssf\prime}}}
\hspace{-1.5mm}d\mu_{\hspace{.2mm}
n^{\ssf\prime}}\ssb(x_{n^{\ssf\prime}}\ssb)
\end{aligned}\end{equation*}
is bounded, uniformly in ${\mathbf{y}}^{\ssf\prime}$.
On the other hand,
let us prove the uniform boundedness of
\begin{equation*}
{\mathbf{I}}^{\ssf\prime\prime}({\mathbf{y}}^{\ssf\prime\prime})\ssf
=\int_{\ssf\R^{n^{\prime\prime}}}\hspace{-1mm}
d{\boldsymbol{\mu}}^{\ssf\prime\prime}(\mathbf{x}^{\ssf\prime\prime})\,
\sup\nolimits_{\,0<t<9\ssf}
t^{\frac{n^\prime}2}\,{\mathbf{H}}_{\,t}^{\ssf\prime\prime}
(\mathbf{x}^{\ssf\prime\prime}\ssb,{\mathbf{y}}^{\ssf\prime\prime})\,,
\end{equation*}
when \ssf$n^{\ssf\prime\prime}\!=\ssb n\ssb-\ssb n^{\ssf\prime}\!>\!0$\ssf.
For this purpose, let us deduce from the Gaussian estimate
\begin{equation*}
{\mathbf{H}}_{\,t}^{\ssf\prime\prime}
(\mathbf{x}^{\ssf\prime\prime}\!,{\mathbf{y}}^{\ssf\prime\prime})
\lesssim\boldsymbol{\mu}^{\ssf\prime\prime}(\ssf\mathbf{B}\ssf
(\ssf{\mathbf{y}}^{\ssf\prime\prime}\ssb,\ssb\sqrt{t\,}))^{-1}\,
e^{-\frac{|\ssf\mathbf{x}^{\ssf\prime\prime}\ssb
-\ssf{\mathbf{y}}^{\ssf\prime\prime}|^2}{c\,t}}
\end{equation*}
that
\begin{equation*}
\sup\nolimits_{\,0<t<9\,}
t^{\frac{n^\prime}2}\,{\mathbf{H}}_{\,t}^{\ssf\prime\prime}
(\mathbf{x}^{\ssf\prime\prime}\ssb,{\mathbf{y}}^{\ssf\prime\prime})
\ssf\lesssim\,
|\ssf\mathbf{x}^{\ssf\prime\prime}\hspace{-1mm}
-\ssb{\mathbf{y}}^{\ssf\prime\prime}|\;
\boldsymbol{\mu}^{\ssf\prime\prime}(\ssf\mathbf{B}
\ssf(\ssf{\mathbf{y}}^{\ssf\prime\prime}\ssb,
|\ssf\mathbf{x}^{\ssf\prime\prime}\hspace{-1mm}
-\ssb{\mathbf{y}}^{\ssf\prime\prime}|\ssf))^{-1}\,
e^{-\frac{|\ssf\mathbf{x}^{\ssf\prime\prime}\ssb
-\ssf{\mathbf{y}}^{\ssf\prime\prime}|^2}{18\,c}}\ssf.
\end{equation*}
Assume first that
\ssf$|\ssf\mathbf{x}^{\ssf\prime\prime}\hspace{-1mm}
-\ssb{\mathbf{y}}^{\ssf\prime\prime}|\ssb\ge\!\sqrt{\ssf t\ssf}$
\ssf with \ssf$0\!<\!t\!<\!9$\ssf.
Then, by using \eqref{ComparisonVolumeBall},
\begin{align*}
t^{\frac{n^\prime}2}\,{\mathbf{H}}_{\,t}^{\ssf\prime\prime}
(\mathbf{x}^{\ssf\prime\prime}\ssb,{\mathbf{y}}^{\ssf\prime\prime})
&\lesssim\ssf
t^{\frac{n^\prime}2}\ssf
\tfrac{\boldsymbol{\mu}^{\ssf\prime\prime}
(\ssf\mathbf{B}\ssf(\ssf{\mathbf{y}}^{\ssf\prime\prime}\ssb,
\ssf|\ssf\mathbf{x}^{\ssf\prime\prime}\ssb
-\ssf{\mathbf{y}}^{\ssf\prime\prime}|\ssf))}
{\boldsymbol{\mu}^{\ssf\prime\prime}
(\ssf\mathbf{B}\ssf
(\ssf{\mathbf{y}}^{\ssf\prime\prime}\ssb,\sqrt{t\ssf}\ssf))}\,
\boldsymbol{\mu}^{\ssf\prime\prime}
(\ssf\mathbf{B}\ssf(\ssf{\mathbf{y}}^{\ssf\prime\prime}\ssb,
|\ssf\mathbf{x}^{\ssf\prime\prime}\hspace{-1mm}
-\ssb{\mathbf{y}}^{\ssf\prime\prime}|\ssf))^{-1}\,
e^{-\frac{|\ssf\mathbf{x}^{\ssf\prime\prime}\ssb
-\ssf{\mathbf{y}}^{\ssf\prime\prime}|^2}{c\,t}}\\
&\lesssim\ssf|\ssf\mathbf{x}^{\ssf\prime\prime}\hspace{-1mm}
-\ssb{\mathbf{y}}^{\ssf\prime\prime}|\ssf
\underbrace{\bigl(\tfrac{|\ssf\mathbf{x}^{\ssf\prime\prime}
\ssb-\ssf{\mathbf{y}}^{\ssf\prime\prime}|}{\sqrt{t\,}}
\bigr)^{\mathbf{N}^{\ssf\prime\prime}}\ssf
e^{-\frac{|\ssf\mathbf{x}^{\ssf\prime\prime}\ssb
-\ssf{\mathbf{y}}^{\ssf\prime\prime}|^2}{2\,c\,t}}}_{\lesssim\,1}
\boldsymbol{\mu}^{\ssf\prime\prime}
(\ssf\mathbf{B}\ssf(\ssf{\mathbf{y}}^{\ssf\prime\prime}\ssb,
|\ssf\mathbf{x}^{\ssf\prime\prime}\hspace{-1mm}
-\ssb{\mathbf{y}}^{\ssf\prime\prime}|\ssf))^{-1}\,
e^{-\frac{|\ssf\mathbf{x}^{\ssf\prime\prime}\ssb
-\ssf{\mathbf{y}}^{\ssf\prime\prime}|^2}{18\,c}}\ssf.
\end{align*}
Assume next that
\ssf$0\ssb<\ssb
|\ssf\mathbf{x}^{\ssf\prime\prime}\hspace{-1mm}
-\ssb{\mathbf{y}}^{\ssf\prime\prime}|
\ssb\le\!\sqrt{\ssf t\ssf}\,(\ssf\le\ssb3\ssf)$\ssf.
Then, by using again \eqref{ComparisonVolumeBall},
\begin{align*}
t^{\frac{n^\prime}2}\,{\mathbf{H}}_{\,t}^{\ssf\prime\prime}
(\mathbf{x}^{\ssf\prime\prime}\ssb,{\mathbf{y}}^{\ssf\prime\prime})
&\lesssim\ssf t^{\frac{n^\prime}2}
\overbrace{\tfrac{\boldsymbol{\mu}^{\ssf\prime\prime}
(\ssf\mathbf{B}\ssf(\ssf{\mathbf{y}}^{\ssf\prime\prime}\ssb,
\ssf|\ssf\mathbf{x}^{\ssf\prime\prime}\ssb
-\ssf{\mathbf{y}}^{\ssf\prime\prime}|\ssf))}
{\boldsymbol{\mu}^{\ssf\prime\prime}(\ssf\mathbf{B}\ssf
(\ssf{\mathbf{y}}^{\ssf\prime\prime}\ssb,\sqrt{t\ssf}\ssf))}
}^{\lesssim\;\bigl(\frac{|\ssf\mathbf{x}^{\ssf\prime\prime}
\!-\ssf{\mathbf{y}}^{\ssf\prime\prime}|}
{\sqrt{t\,}}\bigr)^{\ssb n^{\prime\prime}}}\ssf
\boldsymbol{\mu}^{\ssf\prime\prime}
(\ssf\mathbf{B}\ssf(\ssf{\mathbf{y}}^{\ssf\prime\prime}\ssb,
|\ssf\mathbf{x}^{\ssf\prime\prime}\hspace{-1mm}
-\ssb{\mathbf{y}}^{\ssf\prime\prime}|\ssf))^{-1}
\overbrace{
e^{-\frac{|\ssf\mathbf{x}^{\ssf\prime\prime}\ssb
-\ssf{\mathbf{y}}^{\ssf\prime\prime}|^2}{c\,t}}
}^{\asymp\;1\vphantom{\frac00}}\\
&\lesssim\ssf\underbrace{\vphantom{\tfrac{x^\prime}{\sqrt{t}}}
t^{\frac{n^\prime\!-\ssf1}2}
}_{\lesssim\;1}\underbrace{
\bigl(\tfrac{|\ssf\mathbf{x}^{\ssf\prime\prime}
\ssb-\ssf{\mathbf{y}}^{\ssf\prime\prime}|}
{\sqrt{t\,}}\bigr)^{\ssb n^{\prime\prime}\ssb-\ssf1}
}_{\lesssim\;1}\ssf
|\ssf\mathbf{x}^{\ssf\prime\prime}\hspace{-1mm}
-\ssb{\mathbf{y}}^{\ssf\prime\prime}|\;
\boldsymbol{\mu}^{\ssf\prime\prime}
(\ssf\mathbf{B}\ssf(\ssf{\mathbf{y}}^{\ssf\prime\prime}\ssb,
|\ssf\mathbf{x}^{\ssf\prime\prime}\hspace{-1mm}
-\ssb{\mathbf{y}}^{\ssf\prime\prime}|\ssf))^{-1}
\underbrace{\vphantom{\tfrac{x^\prime}{\sqrt{t}}}
e^{-\frac{|\ssf\mathbf{x}^{\ssf\prime\prime}\ssb
-\ssf{\mathbf{y}}^{\ssf\prime\prime}|^2}{18\,c}}
}_{\asymp\;1}\ssf.
\end{align*}
Now that we have estimated
\ssf$t^{\frac{n^\prime}2}\ssf
{\mathbf{H}}_{\,t}^{\ssf\prime\prime}
(\mathbf{x}^{\ssf\prime\prime}\ssb,
{\mathbf{y}}^{\ssf\prime\prime})$\ssf,
let us split up the integral
\begin{equation*}
{\mathbf{I}}^{\ssf\prime\prime}({\mathbf{y}}^{\ssf\prime\prime})
=\sum\nolimits_{\ssf j\in\Z}\ssf
{\mathbf{I}}_{\ssf j}^{\ssf\prime\prime}
({\mathbf{y}}^{\ssf\prime\prime})
\end{equation*}
according to the decomposition
\,$\R^{n^{\prime\prime}}\hspace{-1mm}\smallsetminus\!\{0\}\ssb
=\bigsqcup_{\,j\in\Z}\underbrace{
\{\ssf\mathbf{x}^{\ssf\prime\prime}\hspace{-1.2mm}\in\!\R^{n^{\prime\prime}}|
\,\,2^{\,j-\frac12}\!\le\ssb
|\ssf\mathbf{x}^{\ssf\prime\prime}\hspace{-1mm}-\ssb\mathbf{y}^{\ssf\prime\prime}|
\ssb<\ssb2^{\,j+\frac12}\ssf\}
}_{\Omega_j}$\,.
Let
\linebreak
\vspace{-4mm}

\noindent
us show that
\begin{equation*}
|\,{\mathbf{I}}_{\ssf j}^{\ssf\prime\prime}
({\mathbf{y}}^{\ssf\prime\prime})\ssf|
\lesssim2^{-|j|}\,.
\end{equation*}
If \,$j\ssb\ge\ssb0$\ssf, we have indeed
\begin{equation*}
{\mathbf{I}}_{\ssf j}^{\ssf\prime\prime}
({\mathbf{y}}^{\ssf\prime\prime})
\lesssim\ssb
\int_{\ssf\Omega_j}\hspace{-1mm}
d{\boldsymbol{\mu}}^{\prime\prime}(\mathbf{x}^{\prime\prime})
\underbrace{
\boldsymbol{\mu}^{\ssf\prime\prime}
(\ssf\mathbf{B}\ssf(\ssf{\mathbf{y}}^{\ssf\prime\prime}\ssb,
|\ssf\mathbf{x}^{\ssf\prime\prime}\hspace{-1mm}
-\ssb{\mathbf{y}}^{\ssf\prime\prime}|\ssf))^{-1}
}_{\asymp\;\boldsymbol{\mu}^{\ssf\prime\prime}
(\ssf\mathbf{B}\ssf(\ssf{\mathbf{y}}^{\ssf\prime\prime}\ssb,
\ssf2^{\ssf j}))^{-1}}
\underbrace{
|\ssf\mathbf{x}^{\ssf\prime\prime}\hspace{-1mm}
-\ssb{\mathbf{y}}^{\ssf\prime\prime}|\,
e^{-\frac{|\ssf\mathbf{x}^{\ssf\prime\prime}\ssb
-\ssf{\mathbf{y}}^{\ssf\prime\prime}|^2}{18\,c}}
}_{\lesssim\;2^{-j}}
\lesssim\ssf\underbrace{\tfrac
{\boldsymbol{\mu}^{\ssf\prime\prime}(\ssf\Omega_j)}
{\boldsymbol{\mu}^{\ssf\prime\prime}(\ssf\mathbf{B}\ssf
(\ssf{\mathbf{y}}^{\ssf\prime\prime}\ssb,\ssf2^{\ssf j}))}
}_{\lesssim\;1}\ssf2^{-j}
\end{equation*}
and, if \,$j\ssb\le\ssb0$\ssf,
\begin{equation*}
{\mathbf{I}}_{\ssf j}^{\ssf\prime\prime}
({\mathbf{y}}^{\ssf\prime\prime})
\lesssim\ssb
\int_{\ssf\Omega_j}\hspace{-1mm}
d{\boldsymbol{\mu}}^{\prime\prime}(\mathbf{x}^{\prime\prime})
\underbrace{
\boldsymbol{\mu}^{\ssf\prime\prime}
(\ssf\mathbf{B}\ssf(\ssf{\mathbf{y}}^{\ssf\prime\prime}\ssb,
|\ssf\mathbf{x}^{\ssf\prime\prime}\hspace{-1mm}
-\ssb{\mathbf{y}}^{\ssf\prime\prime}|\ssf))^{-1}
}_{\asymp\;\boldsymbol{\mu}^{\ssf\prime\prime}
(\ssf\mathbf{B}\ssf
(\ssf{\mathbf{y}}^{\ssf\prime\prime}\ssb,\ssf2^{\ssf j}))^{-1}}
\underbrace{
|\ssf\mathbf{x}^{\ssf\prime\prime}\hspace{-1mm}
-\ssb{\mathbf{y}}^{\ssf\prime\prime}|
}_{\lesssim\;2^{\ssf j}}
\underbrace{
e^{-\frac{|\ssf\mathbf{x}^{\ssf\prime\prime}\ssb
-\ssf{\mathbf{y}}^{\ssf\prime\prime}|^2}{18\,c}}
\vphantom{|}}_{\lesssim\;1}
\lesssim\ssf\underbrace{\tfrac
{\boldsymbol{\mu}^{\ssf\prime\prime}(\ssf\Omega_j)}
{\boldsymbol{\mu}^{\ssf\prime\prime}(\ssf\mathbf{B}\ssf
(\ssf{\mathbf{y}}^{\ssf\prime\prime}\ssb,\ssf2^{\ssf j}))}
}_{\lesssim\;1}\ssf2^{\,j}\ssf.
\end{equation*}
By summing up over \ssf$j\!\in\!\Z$\ssf,
we obtain the uniform boundedness of
\ssf${\mathbf{I}}^{\ssf\prime\prime}
({\mathbf{y}}^{\ssf\prime\prime})$\ssf.
\end{proof}

\section{Proof of Theorem \ref{Theorem1}}
\label{ProofTheorem1}

Theorem \ref{Theorem1} relies on the following result due to Uchiyama \cite{Uchiyama}.

\begin{theorem}\label{TheoremUchiyama}
Assume that a set \,$X$ is equipped with
\begin{itemize}
\item[$\bullet$]
a \ssf{\rm quasi-distance} \ssf$\widetilde{d}$
\ssf i.e.~a distance except that
the triangular inequality is replaced by the weaker condition
\vspace{.5mm}

\centerline{$
\widetilde{d}\ssf(x,y)\le A\,\{\ssf\widetilde{d}\ssf(x,z)+\widetilde{d}\ssf(z,y)\ssf\}
\qquad\forall\;x,y,z\!\in\!X\ssf,
$}\vspace{.5mm}

\item[$\bullet$]
a measure \,$\mu$ whose values on quasi-balls satisfy
\vspace{.5mm}

\centerline{$
\frac rA\le\mu\ssf(\widetilde{B}\ssf(x,r))\le r
\qquad\forall\;x\!\in\!X\ssf,\,\forall\;r\!>\!0\,,
$}\vspace{1mm}

\item[$\bullet$]
a continuous kernel \,$K_r(x,y)\!\ge\!0$ such that,
for every \,$r\!>\!0$ and \,$x,y,y^{\ssf\prime}\hspace{-1mm}\in\!X$,
\vspace{1mm}

\begin{itemize}
\item[$\circ$]
$K_r(x,x)\ge\frac1{A\hspace{.5mm}r}$\,,

\item[$\circ$]
$K_r(x,y)\le r^{-1}\ssf\bigl(\ssf
1\ssb+\ssb\frac{\widetilde{d}\ssf(x,\ssf y)}r
\ssf\bigr)^{\ssb-1-\ssf\delta}$\,,

\item[$\circ$]
$\bigl|\ssf K_r(x,y)\ssb-\ssb K_r(x,y^{\ssf\prime})\bigr|\ssb
\le\ssb r^{-1}\ssf\bigl(\ssf1\ssb
+\ssb\frac{\widetilde{d}\ssf(x,\ssf y)}r\ssf\bigr)^{\ssb-1-\ssf2\ssf\delta}\ssf
\bigl(\ssf\frac{\widetilde{d}\ssf(y,\ssf y^{\ssf\prime})}r\ssf\bigr)^{\ssb\delta}$
\ssf when
\,$\widetilde{d}\ssf(y,y^{\ssf\prime})\!
\le\!\frac{r\ssf+\ssf\widetilde{d}\ssf(x,\ssf y)}{4\ssf A}$\,.
\vspace{.5mm}

\end{itemize}
\end{itemize}
Here \,$A\!\ge\!1$ and \,$\delta\!>\!0$\ssf.
Then the following definitions of the Hardy space \,$H^1(X)$
are equivalent\,{\rm:}
\begin{itemize}
\item[$\bullet$]
{\rm Maximal definition\,:}
$H^1(X)$ \ssf consists of all functions \,$f\!\in\!L^1(X)$ such that
\begin{equation*}
K_*f(x)=\ssf\sup\nolimits_{\,r>0}\,
\Bigl|\ssf{\displaystyle\int_X}d\mu(y)\,K_r(x,y)\,f(y)\,\Bigr|
\end{equation*}
belongs to \ssf$L^1(X)$
and the norm \,$\|f\|_{H^1}$ is comparable to \,$\|K_*f\|_{L^1}$.
\item[$\bullet$]
{\rm Atomic definition\,:}
$H^1(X)$ \ssf consists of all functions \,$f\!\in\!L^1(X)$
which can be written as
\,$f\!=\!\sum_{\ssf\ell}\ssb\lambda_{\ssf\ell}\ssf a_{\ssf\ell}$\ssf,
where the \ssf$a_{\ssf\ell}$\ssb's are atoms \,{\rm(\footnotemark)}
and \,$\sum_{\ssf\ell}\ssb|\lambda_{\ssf\ell}|\!<\!+\infty$\ssf,
and the norm \,$\|f\|_{H^1}$ is comparable to the infimum of
\;$\sum_{\ssf\ell}\ssb|\lambda_{\ssf\ell}|$
\ssf over all such representations.
\end{itemize}
\end{theorem}

\footnotetext{
\,Recall that an atom is a measurable function \,$a\ssb:\ssb X\!\to\ssb\C$
\,such that \ssf$a$ \ssf is supported in a quasi-ball \ssf$\widetilde{B}$\ssf,
\,$\|a\|_{L^\infty}\!\lesssim\ssb\mu\ssf(\widetilde{B})^{-1}$
\ssf and \,$\displaystyle\int_X\!d\mu\,a\ssb=\ssb0$\ssf.
}

Going back to \ssf$X\!=\ssb\R^n$,
equipped with the Euclidean distance
\,$d\ssf(\mathbf{x},\mathbf{y})\ssb=\ssb|\ssf\mathbf{x}\ssb-\ssb\mathbf{y}\ssf|$
\ssf and the measure \eqref{ProductMeasure},
set
\begin{equation*}
\widetilde{d}\ssf(\mathbf{x},\mathbf{y})=\inf\boldsymbol{\mu}\ssf(B)
\qquad\forall\;\mathbf{x},\mathbf{y}\!\in\!\R^n,
\end{equation*}
where the infimum is taken over all closed balls \ssf$B$
\ssf containing \ssf$\mathbf{x}$ \ssf and \ssf$\mathbf{y}$,
and
\begin{equation}\label{RelationHtKr}
K_r(\mathbf{x},\mathbf{y})
=\mathbf{H}_{\ssf t}(\mathbf{x},\mathbf{y})\ssf,
\qquad\forall\;r\!>\!0\ssf,\,\forall\;\mathbf{x},\mathbf{y}\!\in\!\R^n,
\end{equation}
where \,$t\ssb=\ssb t\ssf(\mathbf{x},r)$ is defined by
\ssf$\boldsymbol{\mu}\ssf(B\ssf(\mathbf{x},\ssb\sqrt{t\ssf}\ssf))\ssb=\ssb r$.
In Appendixes B and C,
we check that these data satisfy the assumptions of Uchiyama's Theorem with
\,$\delta\ssb=\ssb\frac1{\mathbf{N}}$\ssf.
Actually the conditions in Theorem \ref{TheoremUchiyama} are obtained up to constants
and they can be achieved by considering suitable multiples of
\ssf$\boldsymbol{\mu}$ \ssf and \ssf$K_r(\mathbf{x},\mathbf{y})$\ssf.
Thus the conclusion of Uchiyama's Theorem hold
for the quasi-distance \ssf $\widetilde{d}$ \ssf
and for the maximal operator \ssf$K_*$\ssf.

On the one hand,
\ssf$d$ \ssf and \ssf$\widetilde{d}$ \ssf define the same Hardy space \ssf$H^1$,
as balls and quasi-balls are comparable.
Let us elaborate.
For every \ssf$\mathbf{x},\mathbf{y}\!\in\!\R^n$ and \ssf$t\!>\!0$\ssf,
we have
\begin{equation*}
|\ssf\mathbf{x}\ssb-\ssb\mathbf{y}\ssf|\le\ssb\sqrt{t\,}
\quad\Longrightarrow\quad
\widetilde{d}\ssf(\mathbf{x},\mathbf{y})\ssb\le r
\quad\Longrightarrow\quad
|\ssf\mathbf{x}\ssb-\ssb\mathbf{y}\ssf|\lesssim\ssb\sqrt{t\,}\,,
\end{equation*}
where \ssf$r\ssb=\ssb\boldsymbol{\mu}\ssf(B\ssf(\mathbf{x},\ssb\sqrt{\ssf t\,}))$.
The first implication is an immediate consequence of the definition of \ssf$\widetilde{d}$
\ssf and the second one is obtained by combining
Lemma \ref{LemmaDistances}.(b) in Appendix B
with \eqref{ComparisonVolumeBall} in Appendix A.
Hence there exists a constant \ssf$c\!>\!0$ \ssf such that
\begin{equation*}
B(\mathbf{x},\ssb\sqrt{\ssf t\,})
\subset\widetilde{B}(\mathbf{x},r)
\subset B(\mathbf{x},c\hspace{.4mm}\sqrt{\ssf t\,})
\end{equation*}
and these sets have comparable measures, according to Appendix A.

On the other hand,
the maximal operators \ssf$K_*$ and \ssf$\mathbf{H}_{\ssf*}$ coincide
and they define the same Hardy space \ssf$H^1$
as the heat maximal operator \ssf$\mathbf{h}_{\ssf*}$\ssf,
according to Propositions
\ref{MaximalOperatorError1D} and \ref{MaximalOperatorErrorProduct}.
Indeed, for every $f\!\in\!L^1(\R^n,d\boldsymbol{\mu})$\ssf, the integrals
\begin{equation*}
\int_{\ssf\R^n}\!d\boldsymbol{\mu}(\mathbf{x})\;\mathbf{h}_{\ssf*}f(\mathbf{x})
\quad\text{and}\quad
\int_{\ssf\R^n}\!d\boldsymbol{\mu}(\mathbf{x})\;\mathbf{H}_{\ssf*}f(\mathbf{x})
\end{equation*}
differ at most by a multiple of \,$\|f\|_{L^1}$,
which is itself controlled by either integral above,
as \,$\mathbf{h}_{\ssf t}(\mathbf{x},\mathbf{y})\,d\boldsymbol{\mu}(\mathbf{y})$
and \ssf$\mathbf{H}_{\ssf t}(\mathbf{x},\mathbf{y})\,d\boldsymbol{\mu}(\mathbf{y})$
are approximations of the identity.

In conclusion,
the atomic Hardy space \ssf$H^1$ associated with Euclidean balls
coincide with the Hardy space \ssf$H^1$ defined by
the heat maximal operator \ssf$\mathbf{h}_{\ssf*}$\ssf.

\hfill$\square$

\section{Proof of Theorem \ref{Theorem2}}

The proof of Theorem \ref{Theorem2} requires
some weighted estimates in Dunkl analysis,
which are well-known in the Euclidean setting
corresponding to \ssf$\mathbf{k}\ssb=\ssb0$\ssf.
Let us first prove a weak analog of the Euclidean estimate
\begin{equation*}
\|\ssf(1\!+\!|\boldsymbol{\xi}|)^\sigma\ssf\widehat{f}(\boldsymbol{\xi})\ssf
\|_{L^1(d\ssf\boldsymbol{\xi})}\ssb\lesssim
\|f\|_{\hspace{.1mm}W_{\ssf2}^{\ssf\sigma\ssf+\ssf n/2\ssf+\ssf\epsilon}}\,.
\end{equation*}

\begin{lemma}\label{Lemma1}
For every \,$\ell\!\in\!\N$ and \,$r\!>\!0$\ssf,
there is a constant \,$C\ssb=\ssb C_{\ell,r}\!>\!0$ such that
\begin{equation*}
\sup\nolimits_{\,\boldsymbol{\xi}\in\R^n}\ssf
(1\!+\!|\boldsymbol{\xi}|)^\ell\,|\ssf\mathcal{F}\ssb f(\boldsymbol{\xi})|
\leq C\,\| f\|_{C^\ell}\,,
\end{equation*}
for every \,$f\!\in\!C^\ell(\R^n)$ with \,$\supp f\!\subset\!B(0,r)$\ssf.
\end{lemma}

\begin{proof}
By using the Riemann-Lebesgue lemma for the Fourier transform \eqref{FourierTransform},
we get
\begin{align*}
\sup\nolimits_{\,\boldsymbol{\xi}\in\R^n}\ssf
(1\!+\!|\boldsymbol{\xi}|)^\ell\,|\ssf\mathcal{F}\ssb f(\boldsymbol{\xi})|\ssf
&\lesssim\,\sup\nolimits_{\,\boldsymbol{\xi}\in\R^n}
\Bigl(1\ssb+\ssb\sum\nolimits_{\ssf j=1}^{\,n}\ssb|\ssf\xi_j|^\ell\ssf\Bigr)\,
|\ssf\mathcal{F}\ssb f(\boldsymbol{\xi})|\\
&\le\,\|f\|_{L^1(d\boldsymbol{\mu})}+\sum\nolimits_{j=1}^{\,n}
\|\ssf D_j^{\ssf\ell}f\ssf\|_{L^1(d\boldsymbol{\mu})}\,.
\end{align*}
The last expression is bounded by \ssf$\| f\|_{C^\ell}$
as, by induction on $\ell$\ssf,
$\operatorname{supp}\ssf( D_j^{\ssf\ell}f)\!\subset\!B(0,r)$
and \ssf$\|\ssf D_j^{\ssf\ell}f\ssf\|_{L^\infty}\!
\lesssim\ssb\|f\|_{C^\ell}\ssf.$
\end{proof}

\begin{corollary}\label{Corollary2}
For every \,$\ell\!\in\!\N$\ssf, \ssf$r\!>\!0$ and \,$\epsilon\!>\!0$\ssf,
there is a constant \,$C\ssb=\ssb C_{\ell,r,\epsilon}\!>\!0$ such that
\begin{equation*}
\|\ssf(1\!+\!|\boldsymbol{\xi}|)^{\ssf\ell\ssf-\ssf\mathbf{N}/2\ssf-\ssf\epsilon}
\ssf\mathcal{F}\ssb f(\boldsymbol{\xi})\ssf
\|_{L^2(d\boldsymbol{\mu}(\boldsymbol{\xi}))}\ssb
\leq C\,\| f\|_{\hspace{.1mm}W_{\ssf2}^{\ssf\ell\ssf+\ssf n/2\ssf+\ssf\epsilon}}\,,
\end{equation*}
for every \,$f\!\in\!W_{\ssf2}^{\ssf\ell\ssf+\ssf n/2\ssf+\ssf\epsilon}(\R^n)$
with \,$\supp f\!\subset\!B(0,r)$.
\end{corollary}

\begin{proof}
This result is deduced from Lemma \ref{Lemma1},
by using on the left hand side the finiteness of the integral
\begin{equation*}
\int_{\ssf\R^n}\!d\boldsymbol{\mu}(\boldsymbol{\xi})\,
(1\!+\!|\boldsymbol{\xi}|)^{-\ssf\mathbf{N}\ssf-\ssf2\ssf\epsilon}
\end{equation*}
and on the right hand side the Euclidean Sobolev embedding theorem.
\end{proof}

\begin{proposition}
For every \,$\sigma\!>\!0$\ssf, $r\!>\!0$ and \,$\epsilon\!>\!0$\ssf,
there is a constant \,$C\ssb=\ssb C_{\sigma,r,\epsilon}\!>\!0$ such that
\begin{equation*}
\|\ssf(1\!+\!|\boldsymbol{\xi}|)^\sigma\ssf\mathcal{F}\ssb f(\boldsymbol{\xi})\ssf
\|_{L^2(d\boldsymbol{\mu}(\boldsymbol{\xi}))}\ssb
\leq C\,\|f\|_{\hspace{.1mm}W_{\ssf2}^{\ssf\sigma\hspace{.1mm}+\ssf\epsilon}}\,,
\end{equation*}
for every \,$f\!\in\!W_{\ssf2}^{\ssf\sigma\hspace{.1mm}+\ssf\epsilon}(\R^n)$
with \,$\supp f\!\subset\!B(0,r)$.
\end{proposition}

\begin{proof}
Let $\chi\!\in\!C_c^\infty(\R^n)$.
Following an argument due to Mauceri-Meda \cite{MauceriMeda},
this result is obtained by interpolation between the $L^2$ estimate
\begin{equation*}
\|\ssf\mathcal{F}(\chi f)\ssf\|_{L^2(d\boldsymbol{\mu})}\ssb
=\const\|\ssf\chi f\ssf\|_{L^2(d\boldsymbol{\mu})}
\lesssim\|f\|_{L^2(d\ssf\mathbf{x})}\ssf,
\end{equation*}
which is deduced from Plancherel's formula,
and the following estimate for \ssf$\ell\!\in\!\N$ \ssf large,
which is deduced from Corollary \ref{Corollary2}\,:
\begin{equation*}
\|\ssf(1\!+\!|\boldsymbol{\xi}|
)^{\ssf\ell\ssf-\ssf\mathbf{N}/2\ssf-\ssf\epsilon^{\ssf\prime}}\ssb
\mathcal{F}(\chi f)(\boldsymbol{\xi})\ssf\|_{L^2(d\boldsymbol{\mu}(\boldsymbol{\xi}))}
\ssb\lesssim\|\ssf\chi f\ssf\|_{\hspace{.1mm}
W_{\ssf2}^{\ssf\ell\ssf+\ssf n/2\ssf+\ssf\epsilon^{\ssf\prime}}}
\ssb\lesssim\|f\|_{\hspace{.1mm}
W_{\ssf2}^{\ssf\ell\ssf+\ssf n/2\ssf+\ssf\epsilon^{\ssf\prime}}}\ssf.
\end{equation*}
\end{proof}

By using the Cauchy-Schwartz inequality,
we deduce eventually the following result.

\begin{corollary}\label{Corollary4}
For every \,$\sigma\!>\!0$\ssf, $r\!>\!0$ and \,$\epsilon\!>\!0$\ssf,
there is a constant \,$C\ssb=\ssb C_{\sigma,r,\epsilon}\!>\!0$ such that
\begin{equation*}
\int_{\ssf\R^n}\!d\boldsymbol{\mu}(\boldsymbol{\xi})\,
(1\!+\!|\boldsymbol{\xi}|)^\sigma\,|\mathcal{F}\ssb f(\boldsymbol{\xi})|
\leq C\,\|f\|_{\hspace{.1mm}
W_{\ssf2}^{\ssf\sigma\hspace{.1mm}+\ssf\mathbf{N}/2\ssf+\ssf\epsilon}}\,,
\end{equation*}
for every \,$f\!\in\!
W_{\ssf2}^{\ssf\sigma\hspace{.1mm}+\ssf\mathbf{N}/2\ssf+\ssf\epsilon}(\R^n)$
with \,$\supp f\!\subset\!B(0,r)$.
\end{corollary}

Let us next prove analogs in the Dunkl setting of the Euclidean estimates
\begin{equation*}
\int_{\ssf\R^n}\!d\ssf\mathbf{x}\,
(1\!+\!|\mathbf{x}|)^\delta\,|f\ssb*\ssb g\ssf(\mathbf{x})|\,
\le\int_{\ssf\R^n}\!d\ssf\mathbf{z}\,
(1\!+\!|\mathbf{z}|)^\delta\,|f(\mathbf{z})|\,
\int_{\ssf\R^n}\!d\ssf\mathbf{y}\,
(1\!+\!|\mathbf{y}|)^\delta\,|g(\mathbf{y})|\,,
\end{equation*}
and
\begin{equation*}
\int_{\ssf\R^n\smallsetminus B(\mathbf{y},\ssf r)}\hspace{-1mm}d\ssf\mathbf{x}\,
|f(\mathbf{x}\!-\!\mathbf{y})|\,\lesssim\,r^{-\delta}\;
\|\ssf(1\!+\!|\mathbf{x}|)^\delta f(\mathbf{x})\ssf\|_{L^1(d\ssf\mathbf{x})}\,.
\end{equation*}
Recall that Dunkl translations
are defined via the Fourier transform \eqref{FourierTransform} by
\begin{equation*}
(\tau_{\ssf\mathbf{y}}f)(\mathbf{x})
=\ssf\mathbf{c}_{\ssf\mathbf{k}}^{-1}\!
\int_{\ssf\R^n}\!d\boldsymbol{\mu}(\boldsymbol{\xi})\,\mathcal{F}f(\mathbf{\boldsymbol{\xi}})\,
\mathbf{E}(\mathbf{x},i\ssf\boldsymbol{\xi})\,
\mathbf{E}(\mathbf{y},i\ssf\boldsymbol{\xi})
\end{equation*}
(see \cite{Roesler2, Trimeche, Roesler4, ThangaveluXu})
and have an explicit integral representation
\begin{equation*}
(\tau_{\ssf\mathbf{y}}f)(\mathbf{x})\ssf
=\ssb\int_{\ssf\R^n}\hspace{-1mm}
d\ssf\boldsymbol{\nu}_{\ssb\mathbf{x},\mathbf{y}\vphantom{|}}(\mathbf{z})\,
f(\mathbf{z})\,,
\end{equation*}
in dimension $1$ (see \cite{Roesler1, ThangaveluXu, AmriAnkerSifi})
and hence in the product case.
Specifically,
\begin{equation*}
d\ssf\boldsymbol{\nu}_{\ssb\mathbf{x},\mathbf{y}\vphantom{|}}(\mathbf{z})
=d\ssf\nu_{x_1,y_1}^{\ssf(1)}\ssb(z_1)\ssf\dots\,
d\ssf\nu_{x_n,y_n}^{\ssf(n)}\ssb(z_n)\,,
\end{equation*}
where
\begin{equation*}
d\ssf\nu_{x_j,y_j}^{\ssf(j)}\ssb(z_j)\,
=\,\begin{cases}
\,\nu_j(x_j,y_j,z_j)\,|\ssf z_j|^{\ssf2\ssf k_j}\ssf dz_j
&\text{if \,}x_j,y_j\!\in\ssb\mathbb{R}^*\ssb,\\
\,d\ssf\delta_{y_j}\ssb(z_j)
&\text{if \,}x_j\!=\ssb0\ssf,\\
\,d\ssf\delta_{x_j}\ssb(z_j)
&\text{if \,}y_j\!=\ssb0\ssf,
\end{cases}
\end{equation*}
and
\begin{equation*}\begin{aligned}
\nu_j(x_j,y_j,z_j)\,
&=\,\tfrac{\Gamma(k_j+\frac12)\vphantom{\big|}}
{\sqrt{\pi\ssf}\,2^{\ssf2\ssf k_j}\ssf\Gamma(k_j)\vphantom{\frac00}}\,
\tfrac{(x_j+\ssf y_j+\ssf z_j)\ssf(-\ssf x_j+\ssf y_j+\ssf z_j)\ssf
(x_j-\ssf y_j+\ssf z_j)\vphantom{\big|}}
{x_j\ssf y_j\ssf z_j\vphantom{\frac00}}\\
&\times\,
\tfrac{\{\ssf(|\ssf x_j|+|\ssf y_j|+|\ssf z_j|)\ssf
(-|\ssf x_j|+|\ssf y_j|+|\ssf z_j|)\ssf(|\ssf x_j|-|\ssf y_j|+|\ssf z_j|)\ssf
(|\ssf x_j|+|\ssf y_j|-|\ssf z_j|)\ssf\}^{\ssf k_j-\ssf1}\vphantom{\big|}}
{|\ssf x_j\ssf y_j\ssf z_j|^{\ssf2\ssf k_j-\ssf1}\vphantom{\frac00}}\\
&\times\,\1_{\,\left[\ssf\left|\vphantom{\frac00}|\ssf x_j|-|\ssf y_j|\right|,
\ssf|\ssf x_j|+|\ssf y_j|\ssf\right]}\ssf(|\ssf z_j|)
\vphantom{\tfrac{\frac00}{\frac00}}\,.
\end{aligned}\end{equation*}
Thus \,$\boldsymbol{\nu}_{\mathbf{x},\mathbf{y}\vphantom{|}}$
\ssf is a signed measure,
which is supported in the product
\begin{equation*}
\boldsymbol{\mathcal{I}}_{\hspace{.15mm}\mathbf{x},\mathbf{y}}
=\ssf\mathcal{I}_{\ssf x_1,\ssf y_1}
\ssb\times\ssf\dots\ssf\times\,
\mathcal{I}_{\ssf x_n,\ssf y_n}
\end{equation*}
of one-dimensional sets
\begin{align*}
\mathcal{I}_{\ssf x_j,\ssf y_j}\ssb
&=\left\{\ssf z_j\!\in\!\R\,|
\left||x_j|\ssb-\ssb|y_j|\vphantom{\tfrac00}\right|\ssb
\le\ssb|z_j|\ssb\leq\ssb|x_j|\ssb+\ssb|y_j|\ssf\vphantom{\tfrac00}\right\}\\
&=\left[-|x_j|\!-\!|y_j|\ssf,-\ssb\left||x_j|\!-\!|y_j|\vphantom{\tfrac00}\right|\ssf\right]
\cup\left[\ssf\left||x_j|\!-\!|y_j|\vphantom{\tfrac00}\right|\ssb,|x_j|\!+\!|y_j|\ssf\right]
\end{align*}
and which is generically given by
\begin{equation*}
d\ssf\boldsymbol{\nu}_{\ssb\mathbf{x},\mathbf{y}\vphantom{|}}(\mathbf{z})
=\underbrace{\nu_1(x_1,y_1,z_1)\ssf\dots\,\nu_n(x_n,y_n,z_n)}_{
\boldsymbol{\nu}(\mathbf{x},\mathbf{y},\mathbf{z})}\,
d\ssf\boldsymbol{\mu}(\mathbf{z})\,.
\end{equation*}
Moreover, it is known (see \cite{Roesler1, ThangaveluXu, AmriAnkerSifi}) that
\begin{equation*}
\sup\nolimits_{\,\mathbf{x},\mathbf{y}\in\R^n}\ssb
|\ssf\boldsymbol{\nu}_{\mathbf{x},\mathbf{y}\vphantom{|}}\ssf|\ssf(\ssf\R^n)
<+\infty\,.
\end{equation*}

\begin{lemma}\label{Lemma5}
For every \,$\delta\ssb\ge\ssb0$\ssf,
$L^1((1\!+\!|\mathbf{x}|)^\delta\ssf d\boldsymbol{\mu}(\mathbf{x}))$
is an algebra with respect to the Dunkl convolution product
\begin{equation*}
f\ssb*\ssb g\ssf(\mathbf{x})
=\!\int_{\ssf\R^n}\hspace{-1mm}d\boldsymbol{\mu}(\mathbf{y})\,
(\tau_{-\mathbf{y}}f)(\mathbf{x})\,g(\mathbf{y})\,\,.
\end{equation*}
\end{lemma}

\begin{proof}
By using the symmetries
\begin{equation*}
\boldsymbol{\nu}(\mathbf{x},-\mathbf{y},\mathbf{z})
=\boldsymbol{\nu}(-\mathbf{z},-\mathbf{y},-\mathbf{x})
=\boldsymbol{\nu}(\mathbf{z},\mathbf{y},\mathbf{x})\,,
\end{equation*}
we have
\begin{equation*}
f\ssb*\ssb g\ssf(\mathbf{x})
=\!\int_{\ssf\R^n}\hspace{-1mm}d\boldsymbol{\mu}(\mathbf{z})\,f(\mathbf{z})
\int_{\ssf\R^n}\hspace{-1mm}d\boldsymbol{\mu}(\mathbf{y})\,g(\mathbf{y})\,
\boldsymbol{\nu}(\mathbf{z},\mathbf{y},\mathbf{x})\,.
\end{equation*}
We conclude by estimating
\begin{equation*}
\int_{\ssf\boldsymbol{\mathcal{I}}_{\hspace{.1mm}\mathbf{z},\mathbf{y}}}
\hspace{-1mm}d\boldsymbol{\mu}(\mathbf{x})\,(1\!+\!|\mathbf{x}|)^\delta\,
|\boldsymbol{\nu}(\mathbf{z},\mathbf{y},\mathbf{x})|
\lesssim(1\!+\!|\mathbf{z}|)^\delta\ssf(1\!+\!|\mathbf{y}|)^\delta\ssf.
\end{equation*}
\end{proof}

\begin{lemma}\label{Lemma6}
For every \,$\delta\!>\!0$\ssf, there is a constant \,$C\!>\!0$ such that,
for every \,$\mathbf{y}\in\R^n$ and \,$r\!>\!0$\ssf,
\begin{equation*}
\int_{\ssf\R^n\smallsetminus\mathcal{O}(\mathbf{y},\ssf r)}\hspace{-1mm}
|(\tau_{-\mathbf{y}}f)(\mathbf{x})|\,d\boldsymbol{\mu}(\mathbf{x})
\le C\,r^{-\delta}\,\|f\|_{L^1((1+|\mathbf{x}|)^\delta d\boldsymbol{\mu}(\mathbf{x}))}\,,
\end{equation*}
\vspace{-3.5mm}

\noindent
where
\begin{equation*}
\mathcal{O}(\mathbf{y},r)
=\{\,\mathbf{x}\!\in\!\R^n\,|\;||x_j|\ssb-\ssb|y_j||\ssb\le\ssb r\hspace{2mm}
\forall\;1\hspace{-.4mm}\le\hspace{-.4mm}j\hspace{-.4mm}\le\hspace{-.4mm}n\,\}
\end{equation*}
\vspace{-4mm}

\noindent
is the orbit of the ball \ssf$B(\mathbf{y},r)$
under the group generated by the reflections \eqref{Reflections}.
\end{lemma}

\begin{proof}
As \,$\R^n\!\smallsetminus\!\mathcal{O}(\mathbf{y},r)$
\ssf is contained in the union of the sets
\begin{equation*}
A_j=\{\,\mathbf{x}\!\in\!\R^n\,|\;||x_j|\ssb-\ssb|y_j||\ssb>n^{-1/2}\,r\,\}
\qquad(\ssf j\ssb=\ssb1\ssf,\ssf\dots,n\ssf)\ssf,
\end{equation*}
we have
\begin{equation*}
\int_{\ssf\R^n\smallsetminus\mathcal{O}(\mathbf{y},r)}\hspace{-1mm}
|(\tau_{-\mathbf{y}}f)(\mathbf{x})|\,d\boldsymbol{\mu}(\mathbf{x})\,
\leq\,\sum\nolimits_{\ssf j=1}^{\,n}
\int_{A_j}\hspace{-1mm}d\boldsymbol{\mu}(\mathbf{x})
\int_{\ssf\boldsymbol{\mathcal{I}}_{\mathbf{x},\mathbf{y}}}
\hspace{-1mm}d\boldsymbol{\mu}(\mathbf{z})
\,|\boldsymbol{\nu}(\mathbf{x},-\mathbf{y},\mathbf{z})|\,|f(\mathbf{z})|\,.
\end{equation*}
As
\begin{equation*}
|\mathbf{z}|\ge|z_j|\ge|\ssf|x_j|\ssb-\ssb|y_j|\ssf|>n^{-1/2}\,r
\end{equation*}
when \ssf$\mathbf{x}\!\in\!A_j$
and \ssf$\mathbf{z}\!\in\!\mathcal{I}_{\ssf\mathbf{x},\mathbf{y}}$\ssf,
the latter expression is bounded above by
\begin{equation*}
n^{\ssf\delta/2}\;r^{-\delta}
\int_{\ssf\R^n}\hspace{-1mm}d\boldsymbol{\mu}(\mathbf{z})\,
|\mathbf{z}|^\delta\,|f(\mathbf{z})|
\int_{\ssf\R^n}\hspace{-1mm}d\boldsymbol{\mu}(\mathbf{x})
\,|\boldsymbol{\nu}(\mathbf{x},-\mathbf{y},\mathbf{z})|\,.
\end{equation*}
We conclude by using the uniform estimate
\begin{equation*}
\int_{\ssf\R^n}\hspace{-1mm}d\boldsymbol{\mu}(\mathbf{x})
\,|\boldsymbol{\nu}(\mathbf{x},-\mathbf{y},\mathbf{z})|
=\ssb\int_{\ssf\R^n}\hspace{-1mm}d\boldsymbol{\mu}(\mathbf{x})
\,|\boldsymbol{\nu}(\mathbf{z},\mathbf{y},\mathbf{x})|
\le C\ssf.
\end{equation*}
\end{proof}

Let us turn to the proof of Theorem \ref{Theorem2},
which consists in estimating
\begin{equation}\label{AtomicEstimateOnRn}
\|\ssf\mathbf{h}_{\ssf*}(\mathcal{T}_{\ssf m\ssf}a)\ssf\|_{L^1(d\boldsymbol{\mu})}
\lesssim M\ssf,
\end{equation}
for every atom \ssf$a$ \ssf in the Hardy space \ssf$H^1$.
By rescaling it suffices to prove \eqref{AtomicEstimateOnRn}
for any atom \ssf$a$ \ssf associated with a unit ball \ssf$B(\mathbf{z},1)$\ssf.
As \ssf$\mathbf{h}_{\ssf*}$ and \ssf$\mathcal{T}_{\ssf m}$
are bounded on \ssf$L^2(\R^n,d\boldsymbol{\mu})$\ssf,
we have
\begin{equation*}\label{AtomicEstimateOnO}
\|\ssf\mathbf{h}_{\ssf*}(\mathcal{T}_{\ssf m\ssf}a)\ssf
\|_{L^1(\mathcal{O}(\mathbf{z},\ssf2),\ssf d\boldsymbol{\mu})}
\ssb\lesssim M\ssf.
\end{equation*}
Thus it remains for us to show that
\begin{equation}\label{AtomicEstimateOnRnMinusO}
\|\ssf\mathbf{h}_{\ssf*}(\mathcal{T}_{\ssf m\ssf}a)\ssf
\|_{L^1(\R^n\smallsetminus\mathcal{O}(\mathbf{z},\ssf2),\ssf d\boldsymbol{\mu})}\ssb
\lesssim M\ssf.
\end{equation}
For this purpose,
let us introduce a dyadic partition of unity on the Dunkl transform side.
More precisely, given a smooth radial function $\psi$ on $\R^n$ such that
\begin{equation*}
\supp\psi\subset
\{\,\xi\!\in\!\R^n\,|\,\tfrac12\ssb\le\ssb|\boldsymbol{\xi}|\ssb\le\ssb2\,\}
\quad\text{and}\quad
\sum\nolimits_{\ssf\ell\in\Z}\psi\ssf(2^{-\ell}\boldsymbol{\xi})^2=1
\quad\forall\;\boldsymbol{\xi}\!\in\!\R^n\!\smallsetminus\!\{0\}\,,
\end{equation*}
let us split up
\begin{equation*}
e^{-\ssf t\ssf|\boldsymbol{\xi}|^2}\ssf m(\boldsymbol{\xi})
=\sum\nolimits_{\ssf\ell\in\Z}
\psi(2^{-\ell}\boldsymbol{\xi})\,e^{-t\ssf|\boldsymbol{\xi}|^2}\ssf
\psi(2^{-\ell}\boldsymbol{\xi})\,m(\boldsymbol{\xi})\,.
\end{equation*}
Set
\vspace{-2mm}
\begin{gather*}
m_{\ssf t,\ell}(\boldsymbol{\xi})
=\overbrace{
\psi(\boldsymbol{\xi})\,e^{-\ssf t\ssf|2^{\ssf\ell}\boldsymbol{\xi}|^2}
}^{\psi_{\ssf t,\ell}(\boldsymbol{\xi})}
\overbrace{\vphantom{e^{|^2}}
\psi(\boldsymbol{\xi})\,m(2^{\ssf\ell}\boldsymbol{\xi})
}^{m_\ell(\boldsymbol{\xi})}\ssf,\\
f_{t,\ell}=\mathcal{F}^{-1}(m_{\ssf t,\ell})
=\underbrace{\mathcal{F}^{-1}(\psi_{\ssf t,\ell})}_{g_{\ssf t,\ell}}\ssf
*\ssf\underbrace{\mathcal{F}^{-1}(m_\ell)}_{\hspace{.1mm}w_\ell}.
\end{gather*}
\vspace{-3mm}

\noindent
Then
\,$e^{-\ssf t\ssf|\boldsymbol\xi|^2}m(\boldsymbol{\xi})\ssb
=\ssb{\displaystyle\sum\nolimits_{\ssf\ell\in\Z}}\ssf
m_{\ssf t,\ell}(2^{-\ell}\boldsymbol{\xi})$\ssf.
Consider the  convolution kernel
\begin{equation*}
F_{t,\ell}(\mathbf{x},\mathbf{y})=
\ssf\tau_{-\mathbf{y}}\,\mathcal{F}^{-1}\bigl\{m_{\ssf t,\ell}(2^{-\ell}\ssf.\ssf)\bigr\}(\mathbf x)
=\ssf2^{\ssf\mathbf{N}\ssf\ell}\ssf
(\tau_{-2^{\ssf\ell}\ssf\mathbf{y}}\ssf f_{t,\ell})(2^{\ssf\ell}\mathbf{x})\ssf.
\end{equation*}
\vspace{-4mm}

\begin{lemma}\label{Lemma9}
\begin{itemize}
\item[(a)]
On the one hand, for every \,$0\!\le\!\delta\!<\!\epsilon$\ssf, we have
\begin{equation*}
\int_{\ssf\R^n\smallsetminus\mathcal{O}(\mathbf{z},\ssf2)}\hspace{-1mm}
d\boldsymbol{\mu}(\mathbf{x})\,
\sup_{\ssf t>0}\ssf|\ssf F_{t,\ell}(\mathbf{x},\mathbf{y})\ssf|\,
\lesssim\,M\,2^{\ssf-\ssf\delta\ssf\ell}\ssf
\qquad\forall\;\ell\!\in\!\Z\ssf,
\;\forall\;\mathbf{z}\hspace{-.75mm}\in\!\R^n,
\;\forall\;\mathbf{y}\!\in\!\mathcal{O}(\mathbf{z},1)\ssf.
\end{equation*}
\item[(b)]
On the other hand,
\begin{equation*}
\int_{\ssf\R^n}\!d\boldsymbol{\mu}(\mathbf{x})\,\sup_{t>0}\ssf
|\ssf F_{t,\ell}(\mathbf{x},\mathbf{y})\ssb
-\ssb F_{t,\ell}(\mathbf{x},\mathbf{y}^{\ssf\prime})\ssf|\,
\lesssim\,M\,2^{\ssf\ell}\ssf|\ssf\mathbf{y}\!-\!\mathbf{y}^{\ssf\prime}|\,
\qquad\forall\;\ell\!\in\!\Z\ssf,\;
\forall\;\mathbf{y},\mathbf{y}^{\ssf\prime}\hspace{-1mm}\in\!\R^n.
\end{equation*}
\end{itemize}
\end{lemma}

\begin{proof}
On the one hand, as
\begin{equation*}
\bigl|\ssf\partial_{\boldsymbol{\xi}}^{\ssf\boldsymbol{\alpha}}
\bigl(\psi(\boldsymbol{\xi})\ssf e^{-t\ssf|\boldsymbol{\xi}|^2}\bigr)
\bigr|\le C_{\boldsymbol{\alpha}}
\qquad\forall\;t\!>\!0\ssf,\;\forall\,\boldsymbol{\xi}\!\in\!\R^n,
\end{equation*}
Lemma \ref{Lemma1} yields the estimate
\begin{equation*}
|\ssf g_{\ssf t,\ell}(\mathbf{x})|\leq C_d\,(1\!+\!|\mathbf{x}|)^{-\ssf d}
\qquad\forall\,\mathbf{x}\!\in\!\R^n,
\end{equation*}
which holds for any \ssf$d\!\in\!\N$
\ssf and which is uniform in \ssf$t\!>\!0$ \ssf and \ssf$\ell\!\in\!\Z$\ssf.
On the other hand, Corollary \ref{Corollary4} yields the estimate
\begin{equation*}
\int_{\ssf\R^n}\!d\boldsymbol{\mu}(\mathbf{x})\,(1\!+\!|\mathbf{x}|)^\delta\,
|w_\ell(\mathbf{x})|\ssf\lesssim M\ssf,
\end{equation*}
which holds uniformly in \ssf$\ell\!\in\!\Z$\ssf.
By resuming the proof of Lemma \ref{Lemma5}, we deduce that
\begin{equation}\label{Equation10}
\int_{\ssf\R^n}\!d\boldsymbol{\mu}(\mathbf{x})\,(1\!+\!|\mathbf{x}|)^\delta\ssf
\sup_{t>0}|f_{t,\ell}(\mathbf{x})|\ssf\lesssim M\ssf.
\end{equation}
We reach our first conclusion by rescaling and by using Lemma \ref{Lemma6}\,:
\begin{align*}
&\int_{\ssf\R^n\smallsetminus\mathcal{O}(\mathbf{z},\ssf2)}
\hspace{-1mm}d\boldsymbol{\mu}(\mathbf{x})\,
\sup_{t>0}\ssf|\ssf F_{t,\ell}(\mathbf{x},\mathbf{y})\ssf|\,
\leq\ssf\int_{\ssf\R^n\smallsetminus\mathcal{O}(\mathbf{y},\ssf1)}
\hspace{-1mm}d\boldsymbol{\mu}(\mathbf{x})\,
\sup_{t>0}\ssf|\ssf F_{t,\ell}(\mathbf{x},\mathbf{y})\ssf|\\
&=\ssf\int_{\ssf\R^n\smallsetminus\mathcal{O}(2^{\ssf\ell}\mathbf{y},\ssf2^{\ssf\ell})}
\hspace{-1mm}d\boldsymbol{\mu}(\mathbf{x})\,
\sup_{t>0}\ssf|(\tau_{-2^{\ssf\ell}\mathbf{y}}\ssf f_{t,\ell})(\mathbf{x})|\,
\lesssim\,M\,2^{\ssf-\ssf\delta\ssf\ell}\ssf.
\end{align*}
Let us turn to the proof of (b).
This time we factorize
\vspace{-1.5mm}
\begin{equation*}
m_{\ssf t,\ell}(\boldsymbol{\xi})
=\ssf\underbrace{\overbrace{
\psi(\boldsymbol{\xi})\,
e^{\ssf|\boldsymbol{\xi}|^2}\ssf
e^{-\ssf t\ssf|2^{\ssf\ell}\boldsymbol{\xi}|^2}
}^{\widetilde{\psi}_{\ssf t,\ell}(\boldsymbol{\xi})}
\overbrace{\vphantom{e^{|^2}}
\psi(\boldsymbol{\xi})\,
m(2^{\ssf\ell}\boldsymbol{\xi})
}^{m_\ell(\boldsymbol{\xi})}
}_{\widetilde{m}_{\ssf t,\ell}(\boldsymbol{\xi})}\,
e^{-|\boldsymbol{\xi}|^2}\ssf,\\
\end{equation*}
\vspace{-3mm}

\noindent
and accordingly
\begin{equation*}
f_{t,\ell}
=\ssf\mathcal{F}^{-1}(m_{\ssf t,\ell})\ssf
=\,\underbrace{
\mathcal{F}^{-1}(\widetilde{m}_{\ssf t,\ell})
}_{\widetilde{f}_{t,\ell}}\,
*\,\underbrace{\vphantom{f_f}
\mathcal{F}^{-1}(e^{-|\boldsymbol{\xi}|^2})
}_{\vphantom{\widetilde{f}}h}\,.
\end{equation*}
\vspace{-2mm}

\noindent
On the one hand, by resuming the proof of \eqref{Equation10}, we get
\begin{equation*}
\int_{\ssf\R^n}\!d\boldsymbol{\mu}(\mathbf{x})\,
\sup_{t>0}|\widetilde{f}_{t,\ell}(\mathbf{x})|\ssf\lesssim M\ssf.
\end{equation*}
On the other hand,
\ssf$\mathbf{h}(\mathbf{x},\mathbf{y})=(\tau_{-\mathbf{y}}\ssf h)(\mathbf{x})$
\ssf is the heat kernel at time \ssf$t\ssb=\ssb1$\ssf, which satisfies
\begin{equation*}
\int_{\ssf\R^n}\!d\boldsymbol{\mu} (\mathbf{x})\,
|\ssf\mathbf{h}(\mathbf{x},\mathbf{y})\ssb
-\ssb\mathbf{h}(\mathbf{x},\mathbf{y}^{\ssf\prime})\ssf|
\lesssim|\ssf\mathbf{y}\ssb-\ssb\mathbf{y}^{\ssf\prime}|
\qquad\forall\;\mathbf{y},\mathbf{y}^{\ssf\prime}\hspace{-1mm}\in\!\R^n,
\end{equation*}
according to next lemma.
After rescaling, we reach our second conclusion\,:
\begin{equation*}
\int_{\ssf\R^n}\!d\boldsymbol{\mu}(\mathbf{x})\,\sup_{\ssf t>0}\ssf
|\ssf F_{t,\ell}(\mathbf{x},\mathbf{y})\ssb
-\ssb F_{t,\ell}(\mathbf{x},\mathbf{y}^{\ssf\prime})\ssf|\ssf
\lesssim\ssf M\,2^{\ssf\ell}\,|\ssf\mathbf{y}\ssb-\ssb\mathbf{y}^{\ssf\prime}|\,.
\end{equation*}
\end{proof}

\begin{lemma}
The following gradient estimate holds for the heat kernel{\rm\;:}
\begin{equation*}
\int_{\ssf\R^n}\!d\boldsymbol{\mu} (\mathbf{x})\,
|\ssf\nabla_{\ssb\mathbf{y}\ssf}\mathbf{h}_{\ssf t}(\mathbf{x},\mathbf{y})\ssf|\,
\lesssim\,t^{-\frac12}
\qquad\forall\;t\!>\!0\ssf,\;\forall\;\mathbf{y}\!\in\!\R^n.
\end{equation*}
\end{lemma}

\begin{proof}
We can reduce to the one-dimensional case
and moreover to \ssf$t\ssb=\!1$ \ssf by rescaling.
It follows from our gradient estimates for the heat kernel in dimension $1$
(see Proposition \ref{PropertiesHeatKernel1D}) that
\vspace{-1mm}
\begin{equation*}
\bigl|\ssf\tfrac\partial{\partial y}\ssf h_1(x,y)\ssf\bigr|
\lesssim\tfrac1{1\ssf+\,|x\ssf y|^k}\,e^{-\frac18\ssf(|x|\ssf-\ssf|y|)^2}\ssf.
\end{equation*}

\noindent$\bullet$
\textit{Case 1}\,:
Assume that \,$|y|\ssb\le\ssb2$\ssf. Then
\ssf$|\ssf\partial_{\ssf y\ssf}h_1(x,y)|\ssb\lesssim\ssb e^{-\ssf x^2\ssb/16}$,
hence
\begin{equation*}
\int_{-\infty}^{+\infty}\hspace{-1mm}dx\,|x|^{\ssf2\ssf k}\,
\bigl|\ssf\tfrac\partial{\partial y}\ssf h_1(x,y)\ssf\bigr|
\,\lesssim\,1\,.
\end{equation*}

\noindent$\bullet$
\textit{Case 2}\,:
Assume that \,$|y|\ssb\ge\ssb2$\ssf. Then
\ssf$|x|/|y|\ssb\le\ssb1\ssb+\ssb\frac12\ssf\bigl||x|\ssb-\ssb|y|\bigr|$\ssf,
hence
\begin{equation*}
|x|^{\ssf2\ssf k}\,\bigl|\ssf\tfrac\partial{\partial y}\ssf h_1(x,y)\ssf\bigr|\ssf
\lesssim\ssf\bigl(\tfrac{|x|\vphantom{|_f}}{|y|\vphantom{|^f}}\bigr)^{\ssb k}\,
e^{-\frac18\ssf(|x|\ssf-\ssf|y|)^2}
\lesssim\ssf\bigl(\ssf1\ssb+\ssb\bigl||x|\ssb-\ssb|y|\bigr|\ssf\bigr)^{\ssb k}\,
e^{-\frac18\ssf(|x|\ssf-\ssf|y|)^2}
\lesssim\ssf e^{-\frac1{16}\ssf(|x|\ssf-\ssf|y|)^2}
\end{equation*}
and
\begin{equation*}
\int_{-\infty}^{+\infty}\hspace{-1mm}dx\,|x|^{\ssf2\ssf k}\,
\bigl|\ssf\tfrac\partial{\partial y}\ssf h_1(x,y)\ssf\bigr|\,
\lesssim\int_{\,0}^{+\infty}\hspace{-1mm}dx\,e^{-\frac1{16}\ssf(x\ssf-\ssf|y|)^2}
\lesssim\int_{-\infty}^{+\infty}\hspace{-1mm}dz\,e^{-\frac1{16}\ssf z^2}
\lesssim\,1\,.
\end{equation*}
\end{proof}

\begin{proof}[End of proof of Theorem \ref{Theorem2}]
Let us split up and estimate
\begin{align*}
|\ssf\mathbf{h}_{\ssf*}(\mathcal{T}_{\ssf m\ssf}a)(\mathbf{x})\ssf|
&\leq\sum\nolimits_{\ssf\ell\ge0}|\ssf\mathbf{h}_{\ssf*}
(\mathcal{T}_{\ssf\psi(2^{-\ell}.\ssf)^2\ssf m\ssf}a)(\mathbf{x})|
+\sum\nolimits_{\ssf\ell<0}|\ssf\mathbf{h}_{\ssf*}
(\mathcal{T}_{\ssf\psi(2^{-\ell}.\ssf)^2\ssf m\ssf}a)(\mathbf{x})|\\
&=\sum\nolimits_{\ssf\ell\ge0}\ssf\sup\nolimits_{\ssf t>0}\,
\Bigl|\,\int_{B(\mathbf{z},\ssf1)}\hspace{-1mm}d\boldsymbol{\mu}(\mathbf{y})\,
F_{t,\ell}(\mathbf{x},\mathbf{y})\,a(\mathbf{y})\,\Bigr|\\
&+\sum\nolimits_{\ssf\ell<0}\ssf\sup\nolimits_{\ssf t>0}\,
\Bigl|\,\int_{B(\mathbf{z},\ssf1)}\hspace{-1mm}d\boldsymbol{\mu}(\mathbf{y})\,
\bigl\{F_{t,\ell}(\mathbf{x},\mathbf{y})\!-\!F_{t,\ell}(\mathbf{x},\mathbf{z})\bigr\}\,
a(\mathbf{y})\,\Bigr|\\
&\leq\sum\nolimits_{\ssf\ell\ge0}\int_{B(\mathbf{z},\ssf1)}
\hspace{-1mm}d\boldsymbol{\mu}(\mathbf{y})\;|a(\mathbf{y})|\,
\sup\nolimits_{\ssf t>0}\,\bigl|\ssf F_{t,\ell}(\mathbf{x},\mathbf{y})\ssf\bigr|\\
&+\sum\nolimits_{\ssf\ell<0}\int_{B(\mathbf{z},\ssf1)}
\hspace{-1mm}d\boldsymbol{\mu}(\mathbf{y})\;|a(\mathbf{y})|\,
\sup\nolimits_{\ssf t>0}\,
\bigl|\ssf F_{t,\ell}(\mathbf{x},\mathbf{y})\ssb
-\ssb F_{t,\ell}(\mathbf{x},\mathbf{z})\ssf\bigr|\,.
\end{align*}
Then \eqref{AtomicEstimateOnRnMinusO} follows from Lemma \ref{Lemma9}.
\end{proof}

\begin{example}
The Riesz transforms 
\,$\mathcal R_j\ssb=\ssb D_j(-\mathbf L)^{-1\slash 2}$
in the Dunkl setting 
correspond to the multipliers
\,$\boldsymbol{\xi}_j/|\boldsymbol{\xi}|$\ssf,
up to a constant.
Hence, by Theorem \ref{Theorem2},
they are bounded operators on  the Hardy space~$H^1$.
\end{example}

\section{Appendixes}

\subsection{Appendix A\,: Measure of balls}
${}$\vspace{1mm}

Recall that \,$k_1,\dots,k_n\ssb\ge\ssb0$
\,and \,$\mathbf{N}\!=\ssb n\ssb+\ssb\sum\nolimits_{\ssf j=1}^{\,n}\ssb2\ssf k_j$\ssf.
On \ssf$\R^n$, equipped with the Euclidean distance,
the product measure
\vskip1.5mm

\centerline{\eqref{ProductMeasure}\hfill$
d\boldsymbol{\mu}(\mathbf{x})
=d\mu_1(x_1)\ssf\dots\,d\mu_n(x_n)
=|x_1|^{\ssf2\ssf k_1}\dots\,|x_n|^{\ssf2\ssf k_n}\,
dx_1\ssf\dots\ssf dx_n
$\hfill}
\vskip1.5mm

\noindent
has the following rescaling properties\,:
\begin{equation}\label{RescalingMeasure}
d\boldsymbol{\mu}(\lambda\ssf\mathbf{x})
=\ssf|\lambda|^{\mathbf{N}}\,d\boldsymbol{\mu}(\ssf\mathbf{x})
\qquad\forall\;\lambda\!\in\!\R^*
\end{equation}
and
\begin{equation}\label{RescalingVolumeBall}
\boldsymbol{\mu}\ssf(\ssf\mathbf{B}\ssf(\lambda\ssf\mathbf{x},|\lambda|\ssf r))
=|\lambda|^{\mathbf{N}}\,
\boldsymbol{\mu}\ssf(\ssf\mathbf{B}\ssf(\ssf\mathbf{x},r))
\qquad\forall\;\mathbf{x}\!\in\!\R^n,\,\forall\;\lambda\!\in\!\R^*.
\end{equation}
Moreover
\begin{equation}\label{VolumeBall}
\boldsymbol{\mu}\ssf(\ssf\mathbf{B}\ssf(\ssf\mathbf{x},r))
\asymp\ssf r^{\ssf n}\ssf\prod\nolimits_{\ssf j=1}^{\,n}\ssf
(\ssf|x_j|\ssb+\ssb r\ssf)^{2\ssf k_j}.
\end{equation}
Hence
\begin{equation}\label{ComparisonVolumeBall}
\bigl(\tfrac Rr\bigr)^{\ssb n}\ssb
\lesssim\tfrac{\boldsymbol{\mu}\ssf(\ssf\mathbf{B}(\mathbf{x},\ssf R\ssf))}
{\boldsymbol{\mu}\ssf(\ssf\mathbf{B}(\mathbf{x},\ssf r))}
\lesssim\ssb\bigl(\tfrac Rr\bigr)^{\ssb\mathbf{N}}
\qquad\forall\;\mathbf{x}\!\in\!\R^n,\,\forall\;R\ssb\ge\ssb r\!>\ssb0\ssf.
\end{equation}
In particular, \ssf$\boldsymbol{\mu}$ \ssf is doubling, i.e.,
\begin{equation}\label{doubling}
\boldsymbol{\mu}\ssf(\ssf\mathbf{B}\ssf(\ssf\mathbf{x},2\ssf r))
\asymp\boldsymbol{\mu}\ssf(\ssf\mathbf{B}\ssf(\ssf\mathbf{x},r))
\qquad\forall\;\mathbf{x}\!\in\!\R^n,\,\forall\;r\!>\ssb0\ssf.
\end{equation}
Let us prove \eqref{VolumeBall} and \eqref{ComparisonVolumeBall}.
In dimension \ssf$n\ssb=\ssb1$\ssf, we have
\begin{equation*}
\mu\ssf(B(x,r))\,
=\int_{\ssf|x|-r}^{\ssf|x|+r}\hspace{-1mm}dy\,|y|^{\ssf2\ssf k}\,.
\end{equation*}
On the one hand, if \,$r\ssb\le\ssb\frac{|x|}2$\ssf, we deduce that
\begin{equation*}
\mu\ssf(B(x,r))\ssf
\asymp\ssf|x|^{\ssf2\ssf k}\int_{\ssf|x|-r}^{\ssf|x|+r}\hspace{-1mm}dy\,
\asymp\ssf|x|^{\ssf2\ssf k}\,r\,.
\end{equation*}
On the other hand, if \,$|x|\ssb\le\ssb2\ssf r$\ssf,
we estimate from above
\begin{equation*}
\mu\ssf(B(x,r))\ssf
\le\ssb\int_{-r}^{\ssf3\ssf r}\hspace{-1mm}dy\,|y|^{\ssf2\ssf k}
\asymp\,r^{\ssf2\ssf k\ssf+1}
\end{equation*}
and from below
\begin{equation*}
\mu\ssf(B(x,r))\ssf
\ge\ssb\int_{\ssf0}^{\ssf r}\!dy\,y^{\ssf2\ssf k}\ssf
\asymp\,r^{\ssf2\ssf k\ssf+1}\ssf.
\end{equation*}
Thus \,$\mu\ssf(B(x,r))\!\asymp\!(\ssf|x|\!+\!r)^{2\ssf k}\ssf r$
\,in all cases and
\begin{equation*}
\mu\ssf(B(x,r))\ssf
\asymp\ssf\Bigl(\tfrac{|x|\ssf+\ssf R}{|x|\ssf+\,r}\Bigr)^{\!2\ssf k}\ssf\tfrac Rr\,
\asymp\,
\begin{cases}
\,\bigl(\frac Rr\bigr)^{\ssb2\ssf k+1}
&\text{if \,}|x|\ssb\le\ssb r\ssf,\\
\,\bigl(\frac R{|x|}\bigr)^{\ssb2\ssf k}\ssf\frac Rr
&\text{if \,}r\ssb\le\ssb|x|\ssb\le\ssb R\,,\\
\hspace{2.5mm}\frac Rr
&\text{if \,}|x|\ssb\ge\ssb R\,.\\
\end{cases}\end{equation*}
The product case follows from the one-dimensional case,
since the ball \ssf$\mathbf{B}\ssf(\ssf\mathbf{x},r)$
\ssf and the cube
\begin{equation*}
\mathbf{Q}\ssf(\ssf\mathbf{x},r)=
\prod\nolimits_{\ssf j=1}^{\,n}\ssb B(x_j,r)
\end{equation*}
have comparable measures. More precisely, we have
\begin{equation*}
\mathbf{Q}\ssf(\ssf\mathbf{x},\tfrac r{\!\sqrt{\ssf n\,}})
\subset\mathbf{B}\ssf(\ssf\mathbf{x},r)
\subset\mathbf{Q}\ssf(\ssf\mathbf{x},r)\ssf,
\end{equation*}
with
\begin{equation*}
\boldsymbol{\mu}\ssf(\mathbf{Q}\ssf(\ssf\mathbf{x},\tfrac r{\!\sqrt{\ssf n\,}}))
\asymp\boldsymbol{\mu}\ssf(\mathbf{Q}\ssf(\ssf\mathbf{x},r))
\asymp\ssf r^{\ssf n}\ssf\prod\nolimits_{\ssf j=1}^{\,n}(\ssf|x_j|\!+\ssb r)^{2\ssf k_j}.
\end{equation*}

\subsection{Appendix B\,: Distances}
${}$\vspace{1mm}

The following result,
which is used in Section \ref{ProofTheorem1},
is certainly known among specialists.
We include nevertheless a proof,
for lack of reference
and for the reader's convenience.

\begin{lemma}\label{LemmaDistances}
Let \ssf$(X,d,\mu)$ be a metric measure space
such that balls have finite positive measure
and satisfy the \ssf{\rm doubling property},
i.e.,
\begin{equation*}
\exists\;C\!>\!0\ssf,
\;\forall\;x\!\in\!X,\;\forall\;r\!>\!0\ssf,
\;\mu\ssf(B(x,2\,r))\le C\,\mu\ssf(B(x,r))\,.
\end{equation*}
Set
\begin{equation*}
\widetilde{d}\ssf(x,y)=\ssf\inf\ssf\mu(B)\ssf,
\end{equation*}
where the infimum is taken over all closed balls \ssf$B$
containing \ssf$x$ and \ssf$y$\ssf.
Then
\begin{itemize}
\item[(a)]
$\,\widetilde{d}$ \ssf is a quasi-distance,
\item[(b)]
$\,\widetilde{d}\ssf(x,y)\ssb\asymp\ssb\mu\ssf(B(x,d\ssf(x,y)))$
\,$\forall\;x,y\!\in\!X$,
\end{itemize}
Moreover,
if the measure \,$\mu$ \ssf has no atoms
and \,$\mu(X)\hspace{-.4mm}=\ssb+\infty$\ssf,
then
\begin{itemize}
\item[(c)]
$\,\mu\ssf(\widetilde{B}(x,r))\ssb\asymp\ssb r$\ssf,
for every \,$x\!\in\!X$ and \,$r\!>\!0$\ssf,
where \,$\widetilde{B}(x,r)$ denotes the closed quasi-ball
with center \ssf$x$ and radius \ssf$r$.
\end{itemize}
\end{lemma}

\begin{proof}
Let us first prove (b).
Set \ssf$R\ssb=\ssb d\ssf(x,y)$\ssf.
On the one hand, we have
\ssf$\widetilde{d}\ssf(x,y)\hspace{-.4mm}\le\ssb\mu\ssf(B(x,R))$\ssf,
as \ssf$x$ \ssf and \ssf$y$ \ssf belong to \ssf$B(x,R)$\ssf.
On the other hand,
if \ssf$x$ \ssf and \ssf$y$ \ssf belong to a ball \ssf$B\ssb=\ssb B(z,r)$\ssf,
then \ssf$R\ssb\le\ssb2\,r$\ssf,
hence \ssf$B(x,R)\ssb\subset\ssb B(z,3\ssf r)$ \ssf and
\ssf$\mu\ssf(B(x,R))\ssb
\le\ssb\mu\ssf(B(z,3\ssf r))\ssb
\asymp\ssb\mu\ssf(B(z,r))$\ssf.
By taking the infimum over all balls \ssf$B$
\ssf containing both \ssf$x$ \ssf and \ssf$y$\ssf,
we conclude that
\ssf$\mu\ssf(B(x,R))\ssb
\lesssim\ssb\widetilde{d}\ssf(x,y)$\ssf.
Let us next prove (a).
For every \ssf$x,y,z\!\in\!X$,
we have \ssf$d\ssf(x,y)\ssb\le\ssb d\ssf(x,z)\ssb+\ssb d\ssf(z,y)$\ssf.
Assume that \ssf$r\ssb
=\ssb d\ssf(x,z)\ssb
\ge\ssb d\ssf(z,y)$\ssf.
Then \ssf$x,y\!\in\!B(z,r)$\ssf.
By using (b), we conclude that
\begin{equation*}
\widetilde{d}\ssf(x,y)
\le\mu\ssf(B(z,r)
\asymp\widetilde{d}\ssf(z,x)
\le\max\,\{\ssf\widetilde{d}\ssf(x,z),\widetilde{d}\ssf(z,y)\}
\le\widetilde{d}\ssf(x,z)\ssb+\ssb\widetilde{d}\ssf(z,y)\,.
\end{equation*}
Let us eventually prove (c).
Given \ssf$x\!\in\!X$, notice that \ssf$\mu\ssf(B(x,r))$
\ssf is an increasing c\`adl\`ag function of \,$r\!\in\!(\ssf0,+\infty)$
\ssf such that
\begin{equation*}\begin{cases}
\;\mu\ssf(B(x,r))\searrow0
&\text{as \;}r\ssb\searrow0\,,\\
\;\mu\ssf(B(x,r))\nearrow+\infty
&\text{as \;}r\ssb\nearrow+\infty\,.\\
\end{cases}\end{equation*}
Here we have used our additional assumptions.
Let \ssf$x\!\in\!X$ and \ssf$r\!>\!0$\ssf.
On the one hand,
for every \ssf$y\!\in\!\widetilde{B}(x,r)$\ssf,
we have \ssf$\mu\ssf(B(x,d\ssf(x,y))\ssb
\asymp\ssb\widetilde{d}\ssf(x,y)\ssb\le\ssb r$\ssf.
Hence
\begin{equation*}
R\ssf=\ssf\sup\,\{\ssf d\ssf(x,y)\ssf|\,y\!\in\!\widetilde{B}(x,r)\ssf\}<+\infty\,.
\end{equation*}
Let \ssf$y\!\in\!\widetilde{B}(x,r)$ \ssf such that
\ssf$d\ssf(x,y)\ssb\ge\ssb\frac R2$\ssf.
Then \ssf$\widetilde{B}(x,r)\ssb
\subset\ssb B(x,R)\ssb\subset\ssb B(x,2\,d\ssf(x,y))$\ssf.
Hence
\begin{equation*}
\mu\ssf(\widetilde{B}(x,r))\le\ssf\mu\ssf(B(x,2\,d\ssf(x,y))
\asymp\ssf\mu\ssf(B(x,d\ssf(x,y))\asymp\ssf\widetilde{d}\ssf(x,y)\le\ssf r\,.
\end{equation*}
On the other hand,
\begin{equation*}
T\ssf=\ssf\inf\,\{\ssf t\!>\!0\,|\,\mu\ssf(B(x,t))\!\ge\!r\ssf\}>0\,.
\end{equation*}
As \ssf$\mu\ssf(B(x,\frac T2))\!<\!r$,
we have \ssf$\widetilde{d}\ssf(x,y)\!<\!r$\ssf,
for every \ssf$y\!\in\!B(x,\frac T2)$\ssf,
hence \ssf$B(x,\frac T2)\!\subset\!\widetilde{B}(x,r)$\ssf.
Consequently,
\begin{equation*}
r\le\ssf\mu\ssf(B(x,T))\asymp\ssf\mu\ssf(B(x,\tfrac T2))
\le\ssf\mu\ssf(\widetilde{B}(x,r))\ssf.
\end{equation*}
\end{proof}

\subsection{Appendix C\,: Kernel bounds}
${}$\vspace{1mm}

Recall from Section \ref{ProofTheorem1} that
the kernels \ssf$K_r(\mathbf{x},\mathbf{y})$ and
\,$\mathbf{H}_{\ssf t}(\mathbf{x},\mathbf{y})$ are related by
\vspace{1.5mm}

\centerline{\eqref{RelationHtKr}\hfill$
K_r(\mathbf{x},\mathbf{y})
=\ssf\mathbf{H}_{\ssf t}(\mathbf{x},\mathbf{y})\,,
$\hfill}\vspace{1.5mm}

\noindent
where \,$r\ssb=\ssb\boldsymbol{\mu}\ssf(B\ssf(\mathbf{x},\ssb\sqrt{t\ssf}\ssf))$\ssf.
In this appendix, we check that the Gaussian estimates of
\ssf$\mathbf{H}_{\ssf t}(\mathbf{x},\mathbf{y})$
in Theorem \ref{EstimatesTruncatedHeatKernelProduct}
imply the estimates of \ssf$K_r(\mathbf{x},\mathbf{y})$
required in Uchiyama's Theorem
(Theorem \ref{TheoremUchiyama}).
This result is certainly well-known among specialists.
We include nevertheless a proof,
for lack of reference and for the reader's convenience.

According to Appendices A and B,
we may consider the quasi-distance \ssf$\widetilde{d}$ \ssf on \ssf$\R^n$
associated with the Euclidean distance
\,$d\ssf(\mathbf{x},\mathbf{y})\ssb=\ssb|\ssf\mathbf{x}\ssb-\ssb\mathbf{y}\ssf|$
and the product measure \eqref{ProductMeasure}.
The on-diagonal lower estimate
\begin{equation}\label{OnDiagonalLowerEstimateKr}
K_r(\mathbf{x},\mathbf{x})\ge\tfrac{C_1}r
\end{equation}
is an immediate consequence of Theorem \ref{EstimatesTruncatedHeatKernelProduct}.(a).
For every \ssf$\delta\!>\!0$\ssf,
the upper estimate
\begin{equation}\label{UpperEstimateKr}
K_r(\mathbf{x},\mathbf{y})\le\tfrac{C_2}r\,\bigl(\ssf1\ssb
+\ssb\tfrac{\widetilde{d}\ssf(\mathbf{x},\ssf\mathbf{y})}r\ssf\bigr)^{\ssb-1-\delta}
\end{equation}
follows from Theorem \ref{EstimatesTruncatedHeatKernelProduct}.(b),
more precisely by combining
\vspace{-1.5mm}
\begin{equation*}
K_r(\mathbf{x},\mathbf{y})\lesssim\ssf
r^{-1}\,e^{-\frac{|\ssf\mathbf{x}\ssf-\ssf\mathbf{y}\ssf|^2}{c\,t}}
\end{equation*}
\vspace{-6.5mm}

\noindent
with
\begin{equation}\label{GaussianYieldsPolynomial}
\bigl(\ssf1\ssb+\ssb
\tfrac{\widetilde{d}\ssf(\mathbf{x},\ssf\mathbf{y})}r\ssf\bigr)^{\ssb1+\delta}
\ssb\le\bigl(\ssf1\ssb+\ssb
\tfrac{\mu\ssf(B(\mathbf{x},\ssf|\ssf\mathbf{x}\ssf-\ssf\mathbf{y}\ssf|))}
{\mu\ssf(B(\mathbf{x},\sqrt{t\ssf}\ssf))}
\ssf\bigr)^{\ssb1+\delta}\ssb\lesssim
\bigl(\ssf1\ssb+\ssb\tfrac{|\ssf\mathbf{x}\ssf-\ssf\mathbf{y}\ssf|}{\sqrt{t\ssf}}
\ssf\bigr)^{\mathbf{N}\ssf(1+\delta)}\ssb
\lesssim\ssf e^{\frac{|\ssf\mathbf{x}\ssf-\ssf\mathbf{y}\ssf|^2}{c\,t}}\ssf.
\end{equation}
The main problem consists in checking the following Lipschitz estimate.

\begin{lemma}\label{LemmaLipschitzEstimateKr}
There exists \ssf$C_3\!>\!0$
\ssf and, for every \ssf$\delta\!>\!0$\ssf,
there exists \ssf$C_4\!>\!0$ \ssf such that
\begin{equation}\label{LipschitzEstimateKr}
\bigl|\ssf K_r(\mathbf{x},\mathbf{y})\ssb
-\ssb K_r(\mathbf{x},\mathbf{y}^{\ssf\prime})\ssf\bigr|
\le\tfrac{C_4}r\ssf\bigl(\ssf
1\ssb+\ssb\tfrac{\widetilde{d}\ssf(\mathbf{x},\ssf\mathbf{y})}r
\ssf\bigr)^{\ssb-1-\ssf\delta}\ssf
\bigl(\ssf\tfrac{\widetilde{d}\ssf(\mathbf{y},\ssf\mathbf{y}^{\ssf\prime})}r
\ssf\bigr)^{\ssb\frac1{\mathbf{N}}}
\end{equation}
if \;$\widetilde{d}\ssf(\mathbf{y},\mathbf{y}^{\ssf\prime})\ssb
\le\ssb C_3\max\,\{\ssf r,\widetilde{d}\ssf(\mathbf{x},\mathbf{y})\ssf\}$\ssf.
\end{lemma}

\begin{proof}
Let us begin with some observations.
First of all,
\eqref{LipschitzEstimateKr} follows from \eqref{UpperEstimateKr},
as long as
\ssf$\widetilde{d}\ssf(\mathbf{y},\mathbf{y}^{\ssf\prime})\ssb\asymp\ssb r$\ssf.
In this case,
we have indeed
\vspace{-1.5mm}
\begin{equation*}
1+\ssf\tfrac{\widetilde{d}\ssf(\mathbf{x},\ssf\mathbf{y})}r\ssf\asymp\ssf
1+\ssf\tfrac{\widetilde{d}\ssf(\mathbf{x},\ssf\mathbf{y}^{\ssf\prime})}r\,.
\end{equation*}
\vspace{-5.5mm}

\noindent
Next, notice that
\vspace{-1mm}
\begin{equation*}\begin{cases}
\,|\ssf\mathbf{x}\ssb-\ssb\mathbf{y}\ssf|\ssb\lesssim\ssb\sqrt{\ssf t\,}
\;\Longleftrightarrow\;
\widetilde{d}\ssf(\mathbf{x},\mathbf{y})\ssb\lesssim r\ssf,\\
\,|\ssf\mathbf{x}\ssb-\ssb\mathbf{y}\ssf|\ssb\gtrsim\ssb\sqrt{\ssf t\,}
\;\Longleftrightarrow\;
\widetilde{d}\ssf(\mathbf{x},\mathbf{y})\ssb\gtrsim r\ssf.
\end{cases}\end{equation*}
This follows indeed from the estimates
\begin{equation*}
\tfrac{\widetilde{d}\ssf(\mathbf{x},\ssf\mathbf{y})}r\asymp
\tfrac{\mu\ssf(B(\mathbf{x},\ssf|\ssf\mathbf{x}\ssf-\ssf\mathbf{y}\ssf|))}
{\mu\ssf(B(\mathbf{x},\sqrt{t\ssf}\ssf))}
\end{equation*}
\vspace{-5.5mm}

\noindent
and
\begin{equation*}
\bigl(\tfrac Rr\bigr)^{\ssb n}\ssb
\lesssim\tfrac{\mu\ssf(B(\mathbf{x},\ssf R))}{\mu\ssf(B(\mathbf{x},\ssf r))}
\lesssim\bigl(\tfrac Rr\bigr)^{\ssb\mathbf{N}}
\quad\text{if \,}r\!\lesssim\!R\ssf.
\end{equation*}
\vspace{-5.5mm}

\noindent
Similarly, we have
\vspace{-1mm}
\begin{equation*}
\,|\ssf\mathbf{y}\ssb-\ssb\mathbf{y}^{\ssf\prime}|\ssb
\lesssim\ssb|\ssf\mathbf{y}\ssb-\ssb\mathbf{x}\ssf|
\;\Longleftrightarrow\;
\widetilde{d}\ssf(\mathbf{y},\mathbf{y}^{\ssf\prime})\ssb
\lesssim\ssb\widetilde{d}\ssf(\mathbf{y},\mathbf{x})\ssf.
\end{equation*}
In particular, there exists \ssf$C_3\!>\!0$ \ssf such that
\vspace{-1mm}
\begin{equation*}
\,|\ssf\mathbf{y}\ssb-\ssb\mathbf{y}^{\ssf\prime}|\ssb
\le\ssb\tfrac12\ssf|\ssf\mathbf{x}\ssb-\ssb\mathbf{y}\ssf|
\quad\text{if}\quad
\widetilde{d}\ssf(\mathbf{y},\mathbf{y}^{\ssf\prime})\ssb
\le C_3\,\widetilde{d}\ssf(\mathbf{x},\mathbf{y})\ssf.
\end{equation*}
Let us turn to the proof of \eqref{LipschitzEstimateKr}
and assume first that
\ssf$\widetilde{d}\ssf(\mathbf{x},\mathbf{y})\hspace{-.5mm}\ge\ssb r$\ssf.
In this case,
\ssf$|\ssf\mathbf{x}\ssb-\ssb\mathbf{y}\ssf|\ssb\gtrsim\!\sqrt{\ssf t\,}$
and \ssf$\widetilde{d}\ssf(\mathbf{y},\mathbf{y}^{\ssf\prime})
\ssb\le\ssb C_3\,\widetilde{d}\ssf(\mathbf{x},\mathbf{y})$\ssf,
hence \ssf$|\ssf\mathbf{y}\ssb-\ssb\mathbf{y}^{\ssf\prime}|
\ssb\le\ssb\frac12\ssf|\ssf\mathbf{x}\ssb-\ssb\mathbf{y}\ssf|$\ssf.
Thus, according to Theorem \ref{EstimatesTruncatedHeatKernelProduct}.(d),
\begin{equation*}
|\ssf K_r(\mathbf{x},\mathbf{y})\ssb
-\ssb K_r(\mathbf{x},\mathbf{y}^{\ssf\prime})\ssf|
=|\ssf\mathbf{H}_{\ssf t}(\mathbf{x},\mathbf{y})\ssb
-\ssb\mathbf{H}_{\ssf t}(\mathbf{x},\mathbf{y}^{\ssf\prime})\ssf|
\end{equation*}
is bounded above by
\vspace{-2mm}
\begin{equation*}
\mu\ssf(B(\mathbf{x},\ssb\sqrt{t\ssf}\ssf))^{-1}\,
e^{-\frac{|\ssf\mathbf{x}\ssf-\ssf\mathbf{y}\ssf|^2}{c\,t}}\,
\tfrac{|\ssf\mathbf{y}\ssf-\ssf\mathbf{y}^{\ssf\prime}|}{\sqrt{\ssf t\,}}\,.
\end{equation*}
After substituting \,$r\ssb=\ssb\mu\ssf(B(\mathbf{x},\ssb\sqrt{t\ssf}\ssf))$
and estimating
\vspace{-1mm}
\begin{equation*}
\bigl(\ssf1\ssb+\ssb\tfrac{\widetilde{d}\ssf(\mathbf{x},\ssf\mathbf{y})}r
\ssf\bigr)^{\ssb1+\delta}\ssb\lesssim\ssf
e^{\frac{|\ssf\mathbf{x}\ssf-\ssf\mathbf{y}\ssf|^2}{2\,c\,t}}
\end{equation*}
as in \eqref{GaussianYieldsPolynomial},
it remains for us to show that
\vspace{-1mm}
\begin{equation}\label{LipschitzAuxiliary}
\tfrac{|\ssf\mathbf{y}\ssf-\ssf\mathbf{y}^{\ssf\prime}|}{\sqrt{\ssf t\,}}
\lesssim\bigl(\ssf
\tfrac{\widetilde{d}\ssf(\mathbf{y},\ssf\mathbf{y}^{\ssf\prime})}r
\ssf\bigr)^{\ssb\frac1{\mathbf{N}}}
\,e^{\frac{|\ssf\mathbf{x}\ssf-\ssf\mathbf{y}\ssf|^2}{2\,c\,t}}\ssf.
\end{equation}
\vspace{-4mm}

\noindent
If \,$|\ssf\mathbf{y}\ssb-\ssb\mathbf{y}^{\ssf\prime}|\ssb\le\!\sqrt{\ssf t\,}$,
then
\vspace{-1.5mm}
\begin{equation*}
\tfrac{\widetilde{d}\ssf(\mathbf{y},\ssf\mathbf{y}^{\ssf\prime})}r
\asymp\tfrac{\mu\ssf(B(\mathbf{y},\ssf|\ssf\mathbf{y}-\mathbf{y}^{\prime}|))}
{\mu\ssf(B(\mathbf{x},\sqrt{t\ssf}\ssf))}
=\tfrac{\mu\ssf(B(\mathbf{y},\ssf|\ssf\mathbf{y}-\mathbf{y}^{\prime}|))}
{\mu\ssf(B(\mathbf{y},\sqrt{t\ssf}\ssf))}\,
\tfrac{\mu\ssf(B(\mathbf{y},\sqrt{t\ssf}\ssf))}
{\mu\ssf(B(\mathbf{x},\sqrt{t\ssf}\ssf))}
\end{equation*}
\vspace{-5.5mm}

\noindent
with
\begin{equation*}
\tfrac{\mu\ssf(B(\mathbf{y},\ssf|\ssf\mathbf{y}-\mathbf{y}^{\prime}|))}
{\mu\ssf(B(\mathbf{y},\sqrt{t\ssf}\ssf))}
\gtrsim\bigl(\ssf
\tfrac{|\ssf\mathbf{y}\ssf-\ssf\mathbf{y}^{\ssf\prime}|}{\sqrt{t\,}}
\ssf\bigr)^{\ssb\mathbf{N}}
\end{equation*}
\vspace{-5mm}

\noindent
and
\vspace{-.5mm}
\begin{equation*}
\tfrac{\mu\ssf(B(\mathbf{y},\sqrt{t\ssf}\ssf))}
{\mu\ssf(B(\mathbf{x},\sqrt{t\ssf}\ssf))}
\ge\tfrac{\mu\ssf(B(\mathbf{y},\sqrt{t\ssf}\ssf))}
{\mu\ssf(B(\mathbf{y},
\ssf|\ssf\mathbf{x}\ssf-\ssf\mathbf{y}\ssf|\ssf+\ssf\sqrt{t\ssf}\ssf))}
\gtrsim\bigl(\ssf\tfrac{\sqrt{t\,}}
{|\ssf\mathbf{x}\ssf-\ssf\mathbf{y}\ssf|\ssf+\ssf\sqrt{t\,}}
\ssf\bigr)^{\ssb\mathbf{N}}\ssb
=\bigl(\ssf1\ssb
+\ssb\tfrac{|\ssf\mathbf{x}\ssf-\ssf\mathbf{y}\ssf|}{\sqrt{t\,}}\ssf
\bigr)^{\ssb-\mathbf{N}}
\ssb\gtrsim\ssf
e^{-\frac{\mathbf{N}}2\ssf\frac{|\ssf\mathbf{x}\ssf-\ssf\mathbf{y}\ssf|^2}{c\,t}}\ssf.
\end{equation*}
If \,$|\ssf\mathbf{y}\ssb-\ssb\mathbf{y}^{\ssf\prime}|\ssb\ge\!\sqrt{t\,}$,
we argue similarly, estimating this time
\vspace{-1.5mm}
\begin{equation*}
\tfrac{\mu\ssf(B(\mathbf{y},\ssf|\ssf\mathbf{y}-\mathbf{y}^{\prime}|))}
{\mu\ssf(B(\mathbf{y},\sqrt{t\ssf}\ssf))}
\gtrsim\bigl(\ssf
\tfrac{|\ssf\mathbf{y}\ssf-\ssf\mathbf{y}^{\ssf\prime}|}{\sqrt{t\,}}
\ssf\bigr)^{\ssb n}
\gtrsim\bigl(\ssf
\tfrac{|\ssf\mathbf{y}\ssf-\ssf\mathbf{y}^{\ssf\prime}|}{\sqrt{t\,}}
\ssf\bigr)^{\ssb\mathbf{N}}\,
\bigl(\ssf
\tfrac{|\ssf\mathbf{x}\ssf-\ssf\mathbf{y}\ssf|}{\sqrt{t\,}}
\ssf\bigr)^{\ssb-(\mathbf{N}-n)}
\gtrsim\bigl(\ssf
\tfrac{|\ssf\mathbf{y}\ssf-\ssf\mathbf{y}^{\ssf\prime}|}{\sqrt{t\,}}
\ssf\bigr)^{\ssb\mathbf{N}}\,
e^{-\frac{\mathbf{N}}4\ssf\frac{|\ssf\mathbf{x}\ssf-\ssf\mathbf{y}\ssf|^2}{c\,t}}
\end{equation*}
\vspace{-5mm}

\noindent
and
\vspace{-1mm}
\begin{equation*}
\tfrac{\mu\ssf(B(\mathbf{y},\sqrt{t\ssf}\ssf))}
{\mu\ssf(B(\mathbf{x},\sqrt{t\ssf}\ssf))}
\gtrsim e^{-\frac{\mathbf{N}}4\ssf
\frac{|\ssf\mathbf{x}\ssf-\ssf\mathbf{y}\ssf|^2}{c\,t}}\ssf.
\end{equation*}
Assume next that
\ssf$\widetilde{d}\ssf(\mathbf{x},\mathbf{y})\hspace{-.4mm}\le\ssb r$\ssf.
Then \ssf$|\ssf\mathbf{x}\ssb-\ssb\mathbf{y}\ssf|\ssb\lesssim\!\sqrt{\ssf t\,}$,
\ssf$\widetilde{d}\ssf(\mathbf{y},\mathbf{y}^{\ssf\prime})\ssb\le\ssb C_3\,r$
\ssf and \eqref{LipschitzEstimateKr} amounts to
\vspace{-1mm}
\begin{equation*}
\bigl|\ssf K_r(\mathbf{x},\mathbf{y})\ssb
-\ssb K_r(\mathbf{x},\mathbf{y}^{\ssf\prime})\ssf\bigr|
\lesssim r^{-1}\ssf\bigl(\ssf
\tfrac{\widetilde{d}\ssf(\mathbf{y},\ssf\mathbf{y}^{\ssf\prime})}r\ssf
\bigr)^{\ssb\frac1{\mathbf{N}}}\ssf.
\end{equation*}
\vspace{-5mm}

\noindent
According to Theorem \ref{EstimatesTruncatedHeatKernelProduct}.(d),
\vspace{-1mm}
\begin{equation*}
|\ssf K_r(\mathbf{x},\mathbf{y})\ssb
-\ssb K_r(\mathbf{x},\mathbf{y}^{\ssf\prime})\ssf|
=|\ssf\mathbf{H}_{\ssf t}(x,\mathbf{y})\ssb
-\ssb\mathbf{H}_{\ssf t}(x,\mathbf{y}^{\ssf\prime})\ssf|
\lesssim \mu\ssf(B(\mathbf{x},\ssb\sqrt{t\ssf}\ssf))^{-1}\,
\tfrac{|\ssf\mathbf{y}\ssf-\ssf\mathbf{y}^{\ssf\prime}|}{\sqrt{t\,}}\,.
\end{equation*}
\vspace{-5mm}

\noindent
As
\vspace{-1mm}
\begin{equation*}
\mu\ssf(B(\mathbf{y},\ssb\sqrt{t\ssf}\ssf))
\asymp\mu\ssf(B(\mathbf{x},\ssb\sqrt{t\ssf}\ssf))=r\,,
\end{equation*}
\vspace{-5mm}

\noindent
we have
\vspace{-.5mm}
\begin{equation*}
\tfrac{\widetilde{d}\ssf(\mathbf{y},\ssf\mathbf{y}^{\ssf\prime})}r\asymp
\tfrac{\mu\ssf(B(\mathbf{y},\ssf|\ssf\mathbf{y}\ssf-\ssf\mathbf{y}^{\ssf\prime}|))}
{\mu\ssf(B(\mathbf{y},\sqrt{t\ssf}\ssf))}\,.
\end{equation*}
\vspace{-3.5mm}

\noindent
As \ssf$\frac{\widetilde{d}\ssf(\mathbf{y},\ssf\mathbf{y}^{\ssf\prime})}r
\hspace{-.4mm}\le\ssb C_3$ \ssf and \ssf$\frac
{\mu\ssf(B(\mathbf{y},\ssf|\ssf\mathbf{y}\ssf-\ssf\mathbf{y}^{\ssf\prime}|))}
{\mu\ssf(B(\mathbf{y},\sqrt{t\ssf}\ssf))}$
\ssf is bounded from below by a power of
\ssf$\frac{|\ssf\mathbf{y}\ssf-\ssf\mathbf{y}^{\ssf\prime}|}{\sqrt{t\,}}$\ssf,
we deduce first that
\ssf$|\ssf\mathbf{y}\ssb-\ssb\mathbf{y}^{\ssf\prime}|\ssb\lesssim\ssb\sqrt{t\ssf}$
\ssf and next that
\vspace{-1.5mm}
\begin{equation*}
\tfrac{\widetilde{d}\ssf(\mathbf{y},\ssf\mathbf{y}^{\ssf\prime})}r
\gtrsim\bigl(\ssf
\tfrac{|\ssf\mathbf{y}\ssf-\ssf\mathbf{y}^{\ssf\prime}|}{\sqrt{t\,}}
\ssf\bigr)^{\ssb\mathbf{N}}\ssf.
\end{equation*}
\vspace{-5mm}

\noindent
This concludes the proof of Lemma \ref{LemmaLipschitzEstimateKr}.
\end{proof}


\begin{thebibliography}{99}

\bibitem{AbramowitzStegun}
M.~Abramowitz, I.A.~Stegun,
\textit{Handbook of mathematical functions\/}
(\textit{with formulas, graphs, and mathematical tables\/}),
Dover Publ. (1972)

\bibitem{AmriAnkerSifi}
B.~Amri, J.-Ph.~Anker, M.~Sifi,
\textit{Three results in Dunkl analysis\/},
Colloq. Math. \textbf{118} (2010), no.~1, 299--312

\bibitem{BurkholderGundySilverstein}
D.L.~Burkholder, R.F.~Gundy, M.L.~Silverstein,
\textit{A maximal function characterisation of the class $H^p$},
Trans. Amer. Math. Soc. \textbf{157} (1971), 137--153

\bibitem{Cherednik}
I.~Cherednik,
\textit{Double affine Hecke algebras\/},
London Math. Soc. Lect. Note Ser. \textbf{319},
Cambridge Univ. Press (2005)

\bibitem{Coifman}
R.R.~Coifman,
\textit{A real variable characterization of $H^p$},
Studia Math. \textbf{51} (1974), 269--274

\bibitem{CoifmanWeiss}
R.R.~Coifman, G.L.~Weiss,
\textit{Extensions of Hardy spaces and their use in analysis\/},
Bull. Amer. Math. Soc. \textbf{83} (1977), 569--615

\bibitem{Dunkl}
C.F.~Dunkl,
\textit{Differential-difference operators associated to reflection groups\/},
Trans. Amer. Math. \textbf{311} (1989), no.~1, 167--183

\bibitem{DziubanskiPreisnerWrobel}
J.~Dziuba\'nski, M.~Preisner, B.~Wr\'obel,
\textit{Multivariate H\"ormander-type multiplier theorem
for the Hankel transform\/},
J. Fourier Anal. Appl. \textbf{19} (2013), no.~2, 417--437

\bibitem{FeffermanStein}
C.~Fefferman, E.M.~Stein,
\textit{$H^p$ \ssb spaces of several variables\/},
Acta Math. \textbf{129} (1972), 137--193

\bibitem{Hormander}
L.~H\"ormander,
\textit{Estimates for translation invariant operators in $L^p$ spaces\/},
Acta Math. \textbf{104} (1960), 93--140

\bibitem{Macdonald}
I.G.~Macdonald,
\textit{Affine Hecke algebras and orthogonal polynomials\/},
Cambridge Tracts Math. \textbf{157},
Cambridge Univ. Press (2003)

\bibitem{MaciasSegovia}
R.A.~Mac\'ias, C.~Segovia,
\textit{A decomposition into atoms of distributions on spaces of homogeneous type\/},
Adv. in Math. \textbf{33} (1979), 271--309

\bibitem{MauceriMeda}
G.~Mauceri, S.~Meda,
\textit{Vector-valued multipliers on stratified groups\/},
Rev. Mat. Iberoamericana \textbf{6} (1990), no. 3--4, 141--154

\bibitem{NIST}
\textit{NIST Digital Library of Mathematical Functions\/},
http\ssf:\ssf/\!/dlmf.nist.gov

\bibitem{Opdam}
E.M.~Opdam,
\textit{Lecture notes on Dunkl operators
for real and complex reflection groups\/},
Math. Soc. Japan Mem. \textbf{8} (2000)

\bibitem{Roesler1}
M.~R\"osler,
\textit{Bessel--type signed hypergroup on \ssf$\mathbb{R}$}\ssf,
in \textit{Probability measures on groups and related structures XI}
(\textit{Oberwolfach, 1994\/}),
H. Heyer \& al. (eds.), World Sci. Publ. (1995), 292--304

\bibitem{Roesler2}
M.~R\"osler,
\textit{Generalized Hermite polynomials and the heat equation for Dunkl operators\/},
Comm. Math. Phys. \textbf{192} (1998), 519--542

\bibitem{Roesler3}
M.~R\"osler,
\textit{Dunkl operators}: \textit{theory and applications\/},
in \textit{Orthogonal polynomials and special functions\/} (\textit{Leuven, 2002\/}),
Lect. Notes Math. \textbf{1817}, Springer-Verlag (2003), 93--135

\bibitem{Roesler4}
M.~R\"osler,
\textit{A positive radial product formula for the Dunkl kernel\/},
Trans. Amer. Math. Soc. \textbf{355} (2003), no. 6, 2413--2438

\bibitem{RoeslerDeJeu}
M.~R\"osler, M.~de Jeu,
\textit{Asymptotic analysis for the Dunkl kernel\/},
J. Approx. Theory \textbf{119} (2002), 110--126

\bibitem{Stein}
E.M.~Stein,
\textit{Harmonic Analysis}
(\textit{Real\ssf-\ssb Variable Methods, Orthogonality, and Oscillatory Integrals\/}),
Princeton Math. Ser. \textbf{43}, Princeton Univ. Press (1993)

\bibitem{SteinWeiss}
E.M.~Stein, G.L.~Weiss,
\textit{On the theory of harmonic functions of several variables I}
(\textit{the theory of $H^p$ spaces\/}),
Acta Math. \textbf{103} (1960), 25--62

\bibitem{ThangaveluXu}
S. Thangavelu and Y. Xu,
\textit{Convolution operator and maximal function for the Dunkl transform\/},
J. Anal. Math. \textbf{97} (2005), 25--56

\bibitem{Trimeche}
K.~Trimeche,
\textit{Paley--Wiener theorems for the Dunkl transform
and Dunkl translation operators},
Integral Transforms Spec. Funct. \textbf{13} (2002), no. 1, 17--38

\bibitem{Uchiyama}
A.~Uchiyama,
\textit{A maximal function characterization of $H^p$
on the space of homogeneous type\/},
Trans. Amer. Math. Soc. \textbf{262} (1980), no. 2, 579--592

\end{thebibliography}
\end{document}